\documentclass[A4paper]{article}

\usepackage{amssymb}
\usepackage{amsfonts}
\usepackage{amsmath}
\usepackage[english]{babel}
\usepackage{bbm}
\usepackage{xcolor}
\usepackage{graphicx}
\usepackage[utf8x]{inputenc}
\usepackage{lmodern}
\usepackage{multicol}
\usepackage{comment}
\usepackage{caption}
\usepackage{subcaption}
\usepackage{hyperref}

\newcommand{\footremember}[2]{%
    \footnote{#2}
    \newcounter{#1}
    \setcounter{#1}{\value{footnote}}%
}
\newcommand{\footrecall}[1]{%
    \footnotemark[\value{#1}]%
}

\newcommand{\dtri}[1]{\textit{Tri}_{#1}}
\newcommand{\dqua}[1]{\textit{Quad}_{#1}}
\newcommand{\dtet}[1]{\textit{Tet}_{#1}}
\newcommand{\dhex}[1]{\textit{Hex}_{#1}}
\newcommand{\dcad}[1]{\textit{Cad}_{#1}}
\newcommand{\param}{\calK}
\newcommand{\rhos}[1]{\varrho_{2, #1}}
\newcommand{\rhov}[1]{\varrho_{3, #1}}
\newcommand{\LTWO}{L^2}
\newcommand{\HONE}{H^1}
\newcommand{\errL}{\calE_{\LTWO}}
\newcommand{\errH}{\calE_{\HONE}}
\newcommand{\delD}{\Delta\text{DOFs}}
\newcommand{\delL}{\Delta\errL}
\newcommand{\delH}{\Delta\errH}
\newcommand{\delT}{\Delta\textit{T}}
\newcommand{\hh}{h}
\newcommand{\Th}{\calT_{\hh}}
\newcommand{\calE}{\mathcal{E}}
\newcommand{\calF}{\mathcal{F}}
\newcommand{\calI}{\mathcal{I}}

\newcommand{\calK}{\mathcal{K}}
\newcommand{\calT}{\mathcal{T}}
\newcommand{\calV}{\mathcal{V}}
\newcommand{\Pol}{E}
\renewcommand{\P}{E}
\newcommand{\F}{F}  
\newcommand{\E}{e}  
\newcommand{\xv}{\mathbf{x}}
\newcommand{\yv}{\mathbf{y}}
\newcommand{\xvP}{\xv_{\P}} 
\newcommand{\mP}{\ABS{\P}}
\newcommand{\mE}{\ABS{\E}}
\newcommand{\hP}{\hh_{\P}}
\newcommand{\hF}{\hh_{\F}}
\newcommand{\ABS}[1]{\vert#1\vert}
\newcommand{\REAL}{\mathbbm{R}}
\newcommand{\HS}[1] {H^{#1}}
\newcommand{\PS}[1]{\mathbbm{P}_{#1}}
\newcommand{\Polynomial}{\mathbb{P}}
\newcommand{\nlen}{\hspace{-0.2mm}}
\newcommand{\snorm}[2]{\vert{#1}\vert_{#2}}
\newcommand{\norm}[2]{\vert\nlen\vert#1\vert\nlen\vert_{#2}}
\newcommand{\sqbra}[1]{\left[{#1}\right]}
\newcommand{\cubra}[1]{\left\{{#1}\right\}}

\begin{document}

\title{Mesh Optimization for the Virtual Element Method: How Small Can an Agglomerated Mesh Become?}

\author{T. Sorgente\footremember{raise}{RAISE Ecosystem, Genova, Italy (tommaso.sorgente@cnr.it, silvia.biasotti@cnr.it, michela.spagnuolo@cnr.it).}, S. Berrone\footremember{disma}{Dipartimento di Scienze Matematiche ``Giuseppe Luigi Lagrange'', Politecnico di Torino, Corso Duca degli Abruzzi 24, 10129 Torino, Italy (stefano.berrone@polito.it, fabio.vicini@polito.it).}, S. Biasotti\footrecall{raise}{}, G. Manzini\footremember{enrico}{Istituto di Matematica Applicata e Tecnologie Informatiche ``Enrico~Magenes'', Consiglio Nazionale delle Ricerche, Via Ferrata 5/A, 27100 Pavia, Italy (marco.manzini@cnr.it).}, M. Spagnuolo\footrecall{raise}{}, F. Vicini\footrecall{disma}{}}

\maketitle

\begin{abstract}
We present an optimization procedure for generic polygonal or polyhedral meshes, tailored for the Virtual Element Method (VEM).
Once the local quality of the mesh elements is analyzed through a quality indicator specific to the VEM, groups of elements are agglomerated to optimize the global mesh quality.
The resulting discretization is significantly lighter: we can remove up to 80$\%$ of the mesh elements, based on a user-set parameter, thus reducing the number of faces, edges, and vertices.
This results in a drastic reduction of the total number of degrees of freedom associated with a discrete problem defined over the mesh with the VEM, in particular, for high-order formulations.
We show how the VEM convergence rate is preserved in the optimized meshes, and the approximation errors are comparable with those obtained with the original ones.
We observe that the optimization has a regularization effect over low-quality meshes, removing the most pathological elements.
This regularization effect is evident in cases where the original meshes cause the VEM to diverge, while the optimized meshes lead to convergence.
We conclude by showing how the optimization of a real CAD model can be used effectively in the simulation of a time-dependent problem.
\end{abstract}

\textbf{Keywords}:
mesh optimization, mesh agglomeration, mesh quality indicators, virtual element method, optimal rates convergence 


\section{Introduction}
\label{sec:introduction}
During the last decades, solving Partial Differential Equations (PDEs) has experienced a substantial surge in its influence on research, design, and production. 
PDEs are indispensable tools for modelling and analyzing phenomena across physics, engineering, biology, and medicine. 
The predominant methods for solving PDEs, such as the Finite Element Method, rely on suitable descriptions of geometrical entities, like the computational domain and its attributes, generally encoded by a mesh.

Despite extensive research and notable achievements, developing techniques for generating meshes with suitable geometrical properties remains an ongoing endeavour. 
Contemporary meshing algorithms typically generate an initial mesh, ideally dominated by well-shaped elements, and optionally execute optimization steps to enhance the geometrical quality of the elements \cite{lo2014finite}. 
Pivotal to these optimizations is the definition of the concept of ``quality'' of a mesh and its elements \cite{sorgente2023survey}.
Although there is some consensus in the literature regarding the notion of quality for triangular/tetrahedral and quadrangular/hexahedral meshes, defining a universal quality indicator for generic polytopal (i.e., polygonal or polyhedral) meshes remains contentious, as evidenced by the myriad attempts \cite{knupp2001algebraic, stimpson2007verdict, chalmeta2013measuring, zunic2004new, huang2020anisotropic, sorgente2021role, sorgente2021indicator, BERRONE2022103770}.

The main reason behind this difficulty is that the concept of ``geometric quality of a cell'' becomes quite vague when the cell is a generic polytope.
Many classical quality indicators become meaningless when the element has a generic number of vertices or faces, or their extension is not straightforward. 
For instance, indicators based on the Jacobian operator can be extended to generic convex polygons, but are not defined for non-convex ones \cite{knupp2001algebraic}.
Most mesh generators and most numerical schemes get around this problem by only allowing convex or star-shaped elements and classifying the others as unlikely configurations.
Alternatively, many mesh optimization strategies are tied solely to the optimization of the elements' size, or the number of incident edges/faces \cite{sorgente2023survey}.

The development of numerical schemes that support generic polytopal elements \cite{BPVEM, Brezzi-Manzini-Marini-Pietra-Russo:2000, BeiraodaVeiga-Lipnikov-Manzini:2014, sukumar2004conforming, di2020hybrid, bertolazzi2004cell} and the great advantages that they provide over classical methods, make it more and more urgent to devise suitable mesh optimization strategies.
In particular, the Virtual Element Method (VEM) was appositely designed to enable computations over any polytopal cell, as it does not require the explicit computation of the basis functions \cite{BBMR}.
To fully exploit its potentialities, we therefore need to be able to generate and optimize polytopal meshes.
Polytopal mesh generation from scratch is generally addressed through Voronoi tessellations \cite{yan2013efficient}, while polytopal mesh optimization is recently emerging as a topic \cite{antonietti2022refinement, antonietti2023agglomeration, sorgente2023mesh} but mainly for two-dimensional domains.

\medskip
In this study, we tackle the problem of optimizing a given planar or volumetric mesh for the VEM, so that it contains the smallest possible number of elements, while still producing accurate results in a VEM numerical simulation.
This is pursued through the use of a polytopal quality indicator \cite{sorgente2021role, sorgente2021indicator} which analyzes the geometrical shape of the elements and indicates groups of elements that can be merged into a single cell without compromising (and sometimes improving) the global mesh quality. 
The algorithm can be seen as a topological optimization method, which modifies the connectivity without changing the vertices position.
To investigate its true potential, we test it on a collection of appositely built non-optimal meshes.
Indeed, if the quality of the original mesh is already high (e.g., if it mainly contains equilateral triangles), there is no real need to optimize it.
Briefly, our algorithm generates optimized meshes that:
\begin{itemize}
    \item contain up to 80$\%$ less elements than the original mesh, consequently reducing the number of degrees of freedom in the discrete problem defined over it;
    \item produce approximation errors comparable or improved to those produced by the original mesh;
\end{itemize}
and these achievements allow us to:
\begin{itemize}
    \item preserve the optimal VEM convergence rates, both in $\LTWO$ and $\HONE$ error norms;
    \item recover the optimal convergence rate in low-quality meshes, by removing pathological elements;
    \item reduce the total computational time, especially in time-dependent simulations.
\end{itemize}

\medskip
The paper is structured as follows.
In Section~\ref{sec:problem}, we introduce the numerical problem to be solved and the VEM discretization used.
We devote Section~\ref{sec:optimization} to the description of the mesh quality indicator and the mesh quality optimization algorithm.
In Section~\ref{sec:datasets} we build a collection of meshes and optimize them with our algorithm.
These meshes will be used for numerical tests in Section~\ref{sec:convergence}, where we report the performance of the VEM over them.
Finally, in Section~\ref{sec:cad}, we present an application of our algorithm to a real CAD model within the context of a time-dependent problem.
To balance completeness and readability, we collect in \ref{app} all the tables omitted from the other sections.

\subsection{Notation and technicalities}
\label{subsec:introduction:notation}
We adhere to the standard definitions and notations of Sobolev spaces, norms, and seminorms, cf.~\cite{Adams-Fournier:2003}. 
Let $k$ be a nonnegative integer.
The Sobolev space $\HS{k}(\omega)$ consists of all square-integrable functions with all square-integrable weak derivatives up to order $k$ that are defined on the open, bounded, connected subset $\omega$ of $\REAL^{d}$, where $d=1,2,3$.
In the case where $k=0$, we use the notation $\LTWO(\omega)$.
The norm and seminorm in $\HS{k}(\omega)$ are represented by $\norm{\cdot}{k,\omega}$ and $\snorm{\cdot}{k,\omega}$ respectively, and for the inner product in $\LTWO(\omega)$ we use the $(\cdot, \cdot)_{\omega}$ notation.

\begin{table}[!ht]
    \centering
    \caption{Notation summary.} \label{tab:notation}
    \includegraphics[width=0.6\textwidth]{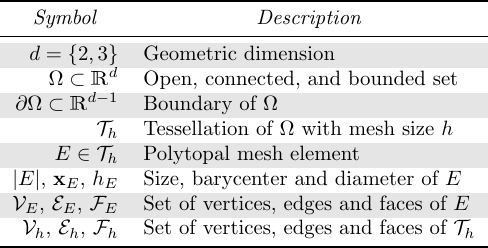}
\end{table}

Table~\ref{tab:notation} reports a brief summary of the notation used in this paper.
We denote by $\Omega \subset \REAL^d$, $d \in \{2, 3\}$ an open, connected, and bounded set, with boundary $\partial \Omega$.
We employ the symbol $\Th$ to represent a tessellation of $\Omega$ with mesh size $\hh$, i.e., a partition of $\Omega$ into non-overlapping polytopal elements $\P$ such that $\bigcup_{\P \in \calT_{\Omega}} \P = \bar{\Omega}$.
A polytopal element $\P$ (a polygon for $d=2$ or a polyhedron for $d=3$) is a compact subset of $\REAL^d$ with boundary $\partial\P$, size (area or volume) $\mP$, barycenter $\xvP$, and diameter $\hP=\sup_{\xv,\yv\in\P}\vert\xv-\yv\vert$.
The mesh size labeling each mesh $\Th$ is defined by $\hh=\max_{\P\in\Th}\hP$.
We denote by $\calV_{\P}$, $\calE_{\P}$, and $\calF_{\P}$ (if $d = 3$) the set of vertices, edges and faces of the element $\P$, respectively.
With the same intention, we use $\calV_{\hh}$, $\calE_{\hh}$, and $\calF_{\hh}$ to indicate the set of vertices, edges and faces (if $d = 3$) of the tesselation $\Th$.

\section{The Virtual Element Discretization}
\label{sec:problem}
In this section, we consider an elliptic second-order differential problem in $\Omega \subset \REAL^d$, $d \in \{2, 3\}$, written in its variational formulation: given the functional space $V$, we look for $u \in V$ the solution of
\begin{equation}
    \label{eq:prob:var}
    a(u, v) = F(v) \quad \forall v \in V,
\end{equation}
where 
\begin{itemize}
\item $V$ is the Sobolev space $H^1(\Omega)$, or $H^1_{\Gamma}(\Omega)$, $\Gamma$ being a non-empty subset of $\partial\Omega$ with a non-zero $(d-1)$ dimensional Lebesgue measure where some set of essential boundary conditions are possibly imposed;
\item $a : V \times V \to \REAL$ is a continuous and coercive bilinear form; 
\item $F : V \to \REAL$ is a continuous linear functional.
\end{itemize}
Under these hypothesis, the Lax-Milgram lemma implies that Problem~\eqref{eq:prob:var} is well-posed; see, e.g., \cite{Brenner-Scott:2008}.

We approach the discretization of such a variational problem through the Virtual Element Method in a $d$-dimensional setting, as described in \cite{AHMAD2013376, VEMGSO}.
Let $\P$ be a generic $d$-dimensional polytopal element of the domain partition $\Th$.
We fix a positive integer $k$ as the order of the discrete approximation.
We denote $\Polynomial^d_k (\P)$ the set of $d$-dimensional polynomials of degree less than or equal to $k$ defined on $\Pol$.
We conventionally assume that $\PS{-1}(\omega)=\{0\}$.
In our implementation of the VEM for numerical tests, we use the standard basis of scaled monomials for all polynomial spaces \cite{BPVEM}.

We introduce two elliptic projection operators on $\Pol$: namely, $\Pi_k^{\nabla} : H^1(\Pol) \to \Polynomial^d_k (\Pol)$ defined as
\begin{equation*}
    \begin{cases}
        \displaystyle
    \left( \nabla p, \nabla(\Pi_k^{\nabla} v - v) \right)_{\Pol} = 0, & \forall p \in \Polynomial^d_k(\Pol),
    \\[0.5em]
    \displaystyle
    \left(1, \Pi_k^{\nabla} v - v  \right)_{\partial \Pol} = 0, &
    \end{cases}
\end{equation*}
and $\Pi_k^0 : L^2(\Pol) \to \Polynomial_k (\Pol)$ as
\begin{equation*}
    \left(p, \Pi_k^{0} v - v\right)_{\Pol} = 0, \quad \forall p \in \Polynomial^d_k(\Pol).
\end{equation*}
On every element $\P$, we consider the local virtual element space of order $k$ that is defined as
\begin{itemize}
    \item if $d = 2$,
    \begin{eqnarray}
        \label{eq:vem:space:f}
        V^{\P}_k &=& \big\{ v\! \in\! H^1(\P) : \Delta v \in \Polynomial^2_k(\P),\ v_{|_{{\E}}}\! \in\! \Polynomial^1_k(\E)\ \forall \E \in \calE_{\P},\ v_{|_{{\partial \P}}} \in C^0(\partial \P),\nonumber\\[0.25em]
        &&\phantom{\big\{ v\! \in\! H^1(\P) : }
        \left(p, v - \Pi_k^{\nabla}v \right)_{\P} = 0\ \forall p \in \Polynomial^2_k(\P)\ /\ \Polynomial^2_{k - 2}(\P) \big\},
    \end{eqnarray}
    \item if $d = 3$,
    \begin{eqnarray}
        \label{eq:vem:space:p}
        V^{\P}_k &=& \big\{ v\! \in\! H^1(\P) : \Delta v \in \Polynomial^3_k(\P),\ v_{|_{{\F}}}\! \in\! V^{\F}_k \ \forall \F \in \calF_{\P},\ v_{|_{{\partial \P}}} \in C^0(\partial \P),\nonumber\\[0.25em]
        &&\phantom{\big\{ v\! \in\! H^1(\P) :}
        \left(p, v - \Pi_k^{\nabla}v \right)_{\P} = 0\ \forall p \in \Polynomial^3_k(\P)\ /\ \Polynomial^3_{k - 2}(\P) \big\}.
    \end{eqnarray}
\end{itemize}
The symbol $\Polynomial^d_k(\Pol)\ /\ \Polynomial^d_{k - 2}(\Pol)$ denotes the set of polynomials in $\Polynomial^d_k(\Pol)$ that are $L^2$-orthogonal to $\Polynomial^d_{k-2}(\Pol)$, and $V^{\F}_k$ is the virtual element space defined on the two-dimensional face $\F$ defined for $d=2$.

With the local spaces of Equations~\eqref{eq:vem:space:f}-\eqref{eq:vem:space:p}, we define the global virtual element space
\begin{equation}
    \label{eq:vem:space}
    V_k = \left\{ v \in C^0(\bar{\Omega}) \cap V : v_{|_{{\P}}} \in V^{\P}_k \ \forall \P \in \Th \right\}.
\end{equation}
Each function $v \in V_k \subset V$ can be uniquely identified, see \cite{BPVEM}, by the set of Degrees of Freedoms (DOFs)
\begin{itemize}
    \item the value of $v$ in each vertex of the set $\calV_{\hh}$;
    \item if $k > 1$, the value of $v$ in the internal $k-1$ Gauss quadrature nodes of each edge of set $\calE_{\hh}$;
    \item if $k > 1$ and $d = 3$, the value of the internal scaled moments of $v$ of order $k - 2$ on each face of set $\calF_{\hh}$;
    \item if $k > 1$, the value of the internal scaled moments of $v$ of order $k - 2$ on each $\P \in \Th$.
\end{itemize}
Both projection operators $\Pi_k^{\nabla}$ and $\Pi_k^{0}$ are computable on every $\P$ from these degrees of freedom, although the construction of $\Pi_k^{\nabla}$ on $\P$ recursively requires the construction of $\Pi_k^{\nabla}$ on every face $\F\subset\partial\P$.

We discretize Problem~\eqref{eq:prob:var} using the local bilinear discrete form $a^{\P}_h : V^{\P}_k \times V^{\P}_k \to \REAL$:
\begin{equation*}
    \label{eq:vem:a}
    a^{\P}_{h}(u, v) := a^{\P}(\Pi_k^{\nabla}u, \Pi_k^{\nabla}v) + S^{\P}(u - \Pi^{\nabla}_k u, v - \Pi^{\nabla}_k v).
\end{equation*}
where 
\begin{itemize}
\item 
$a^{\P}_{h}(\cdot,\cdot)$ is an approximation of the bilinear form $a^{\P}(\cdot, \cdot)$, the restriction of the bilinear form $a(\cdot,\cdot)$ to the element $\P$;
\item the bilinear form $S^{\P}(\cdot, \cdot)$ is the \emph{stabilization term}, which 
can be any computable, symmetric, positive definite bilinear form satisfying the stability condition
\begin{equation*}
    c_{*} a(v, v) \leq S(v, v) \leq c^{*} a(v, v),\ \forall v \in V^{\P}_k \cap \text{ker}(\Pi^{\nabla}_k), 
\end{equation*}
for some pair of real, positive constants $c_*$ and $c^*$.
\end{itemize}
In the implementation for the numerical results, we apply the standard \emph{dofi-dofi} stabilization proposed in \cite{BPVEM}.

This leads to the discrete counterpart of Problem~\eqref{eq:prob:var}: 
\emph{Find $u_h \in V_k$ such that:}
\begin{equation}
    \label{eq:prob:vem}
    a_h(u_h,v_h) := 
    \sum_{\P \in Th} a^{\P}_h(u_h, v_h) = \sum_{\P \in Th} F^{\P}_h(v_h) 
    =: F_h(v_h)
    \quad \forall v_h \in V_k,
\end{equation} 
where, for every $v_h\in V^\P_k$, we set $F^{\P}_h(v_h)=F^\P(\Pi^0_k v_h)$ as a local approximation of $F(v_h)$. 
Other details about $a(\cdot,\cdot)$, $F(\cdot)$, and their virtual element approximation are given in Section~\ref{sec:convergence}.

\section{Mesh Quality Optimization Algorithm}
\label{sec:optimization}
\begin{figure}[!h]
  \centering
  \begin{subfigure}[]{.24\linewidth}
      \centering
      \includegraphics[width=\linewidth]{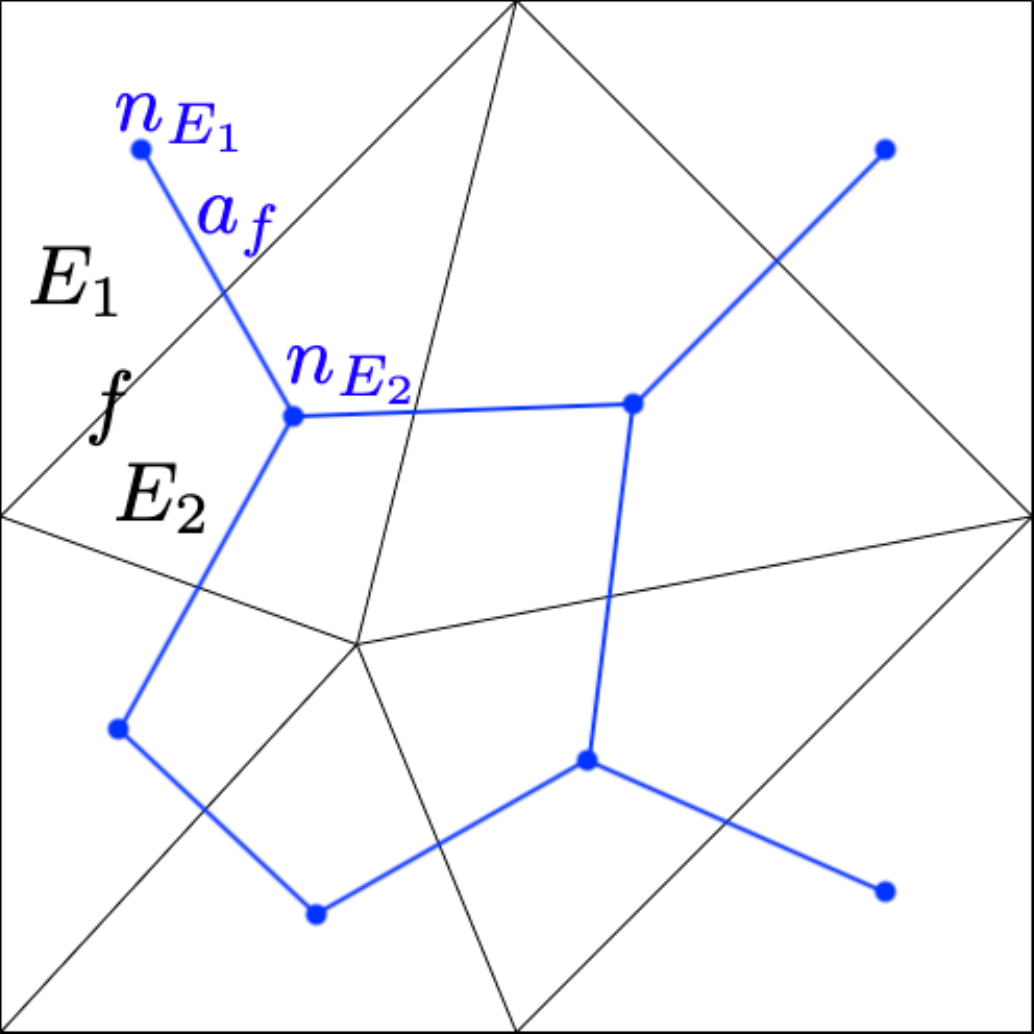} \caption{} \label{fig:pipeline:1}
  \end{subfigure}
  \begin{subfigure}[]{.24\linewidth}
      \centering
      \includegraphics[width=\linewidth]{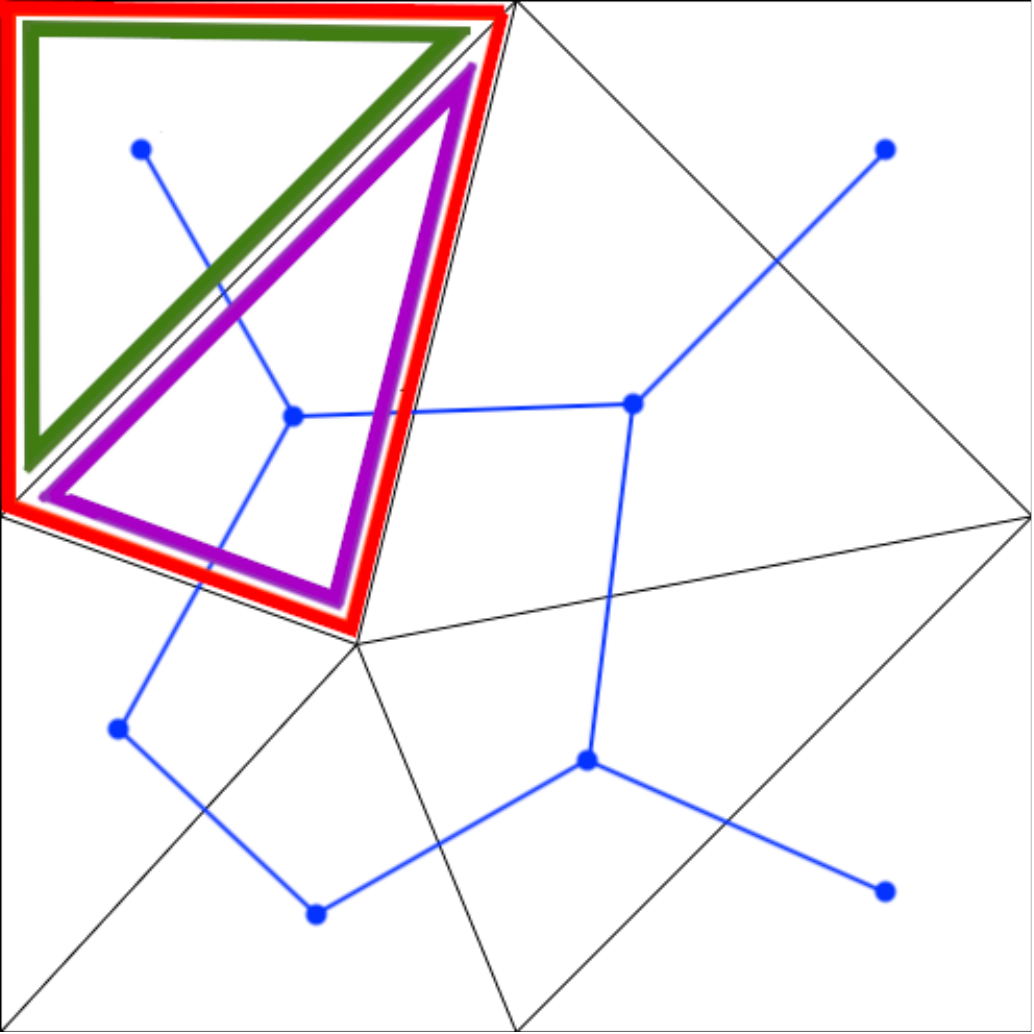} \caption{} \label{fig:pipeline:2}
  \end{subfigure}
  \begin{subfigure}[]{.24\linewidth}
      \centering
      \includegraphics[width=\linewidth]{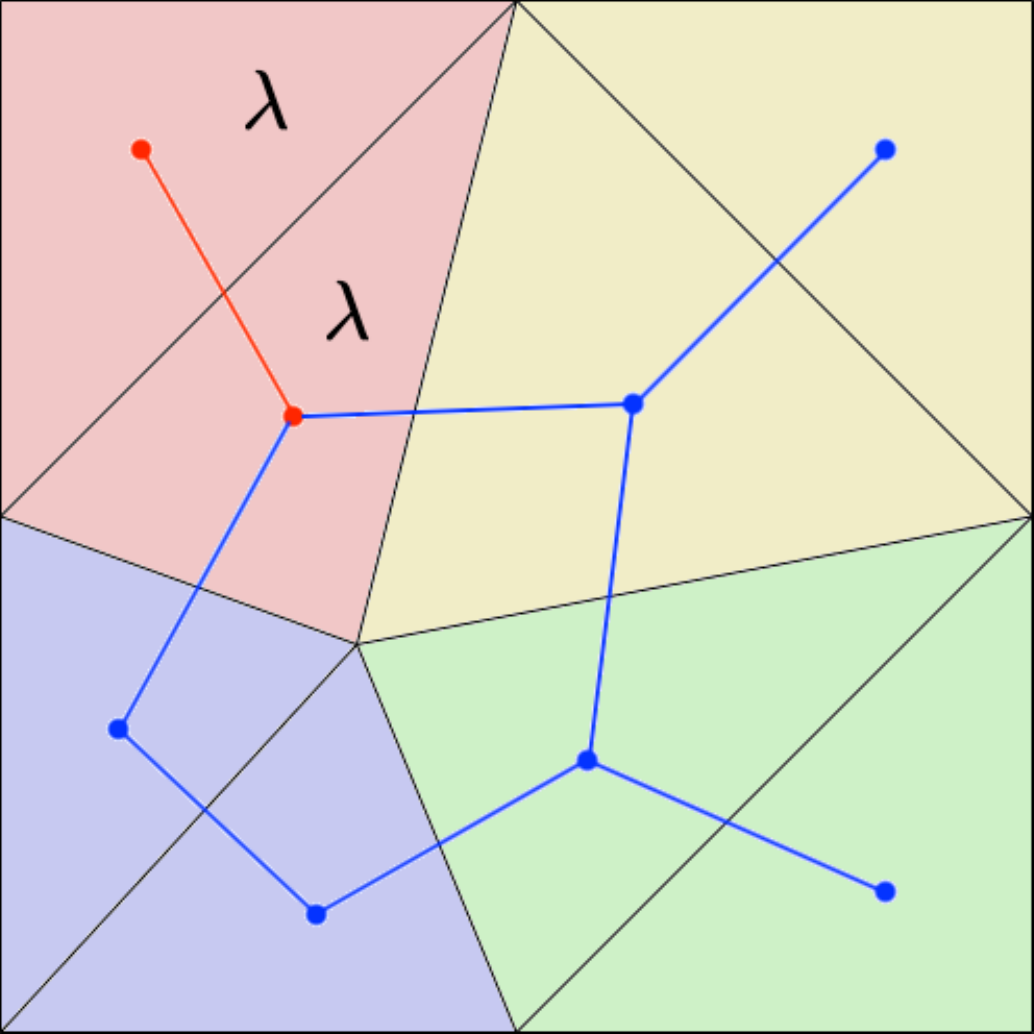} \caption{} \label{fig:pipeline:3}
  \end{subfigure}
  \begin{subfigure}[]{.24\linewidth}
      \centering
      \includegraphics[width=\linewidth]{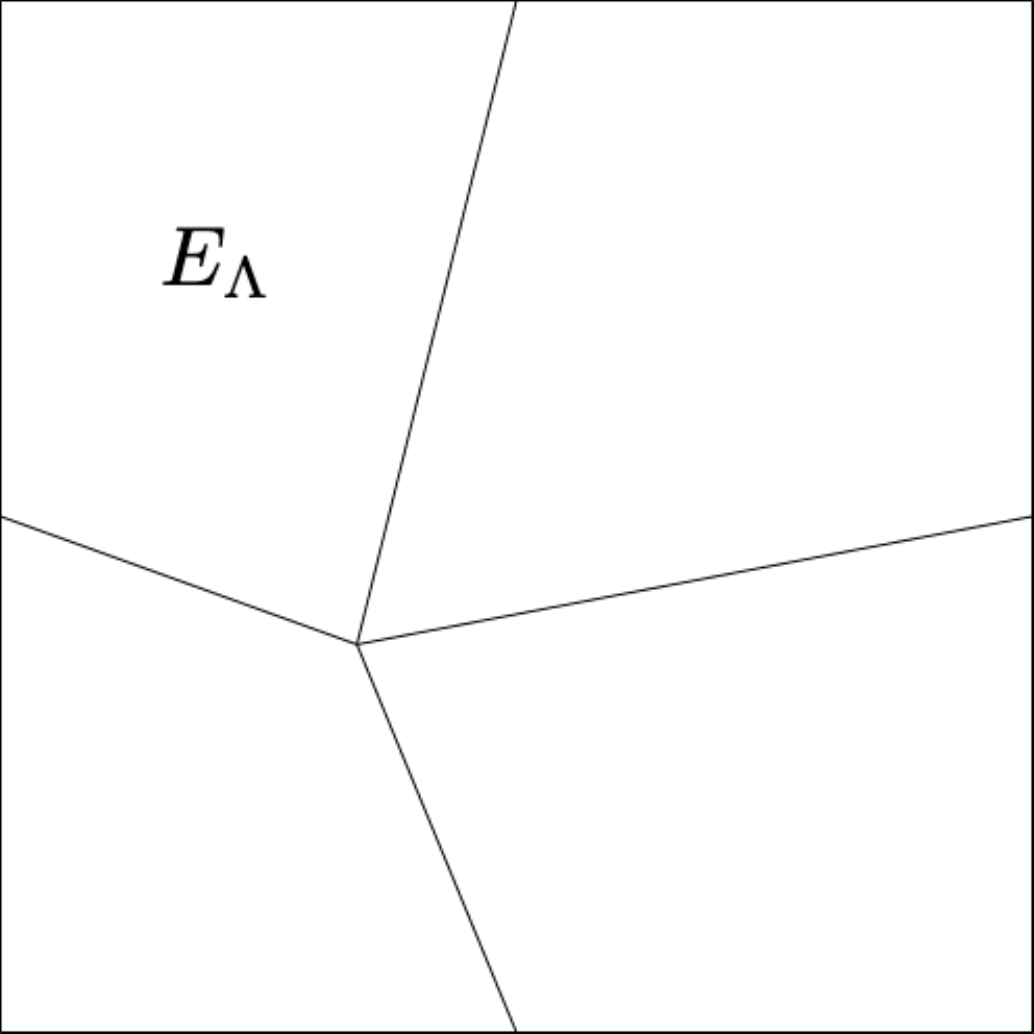} \caption{} \label{fig:pipeline:4}
  \end{subfigure}
  \caption{Optimization pipeline: 
  (a)~input mesh (black), with elements $\P_1, \P_2$ and their shared edge $f$, and mesh dual graph (blue), with nodes $n_{\P_1}, n_{\P_2}$ and arc $a_f$; 
  (b)~weights assignment: $w(n_{\P_1})$ is the quality of the green polygon, $w(n_{\P_2})$ is the quality of the purple polygon and $w(a_f)$ is the quality of the red polygon; 
  (c)~input mesh colored with respect to METIS labeling: nodes $n_{\P_1}$ and $n_{\P_2}$ are assigned the same label $\lambda$; 
  (d)~agglomeration of the elements with the same label: elements $\P_1$ and $\P_2$ become a single element $\P_\Lambda$.}
  \label{fig:pipeline}
\end{figure}

The mesh optimization process involves two distinct steps.
Initially, we evaluate the quality of the mesh elements by employing a local quality indicator.
Subsequently, these elements are categorized and tagged, followed by the agglomeration of elements that bear identical tags.
This methodology is applicable to both polygonal and polyhedral meshes.
We detail the process specifically for three-dimensional meshes ($d=3$), noting any variation that is relevant to two-dimensional cases (d = 2).
The entire workflow, shown in Figure~\ref{fig:pipeline}, has been impleented using the \textit{cinolib} library \cite{livesu2019cinolib}.

\subsection{Quality Indicator}
\label{subsec:optimization:quality}
We borrow the notion of quality from Sorgente et al
\cite{sorgente2021role, sorgente2021indicator}, who propose a mesh
quality indicator specifically tailored for the Virtual Element
Method.
We briefly report the indicator definitions and refer the reader to
the original papers for a more complete discussion.
Given a polygonal mesh $\Th$, for each element
$\P\in\Th$ we have the following indicators:
\begin{align}
\rhos{1}(\P) &= \frac{\lvert kernel(\P)\rvert}{\mP} = 
    \begin{cases} 
    1 & \mbox{if $\P$ is convex} \\
    \in (0,1) & \mbox{if $\P$ is concave and star-shaped} \\
    0 & \mbox{if $\P$ is not star-shaped} \\
    \end{cases}\notag\\
\rhos{2}(\P) &= \frac{\min(\sqrt{\mP},\ \min_{\E\in\calE_{\P}}\mE)} {\hP} \notag\\
\rhos{3}(\P) &= \frac{3}{\#\left\{\E\in\calE_{\P} \right\}} = 
    \begin{cases} 
    1 & \mbox{if $\P$ is a triangle} \\
    \in (0,1) & \mbox{otherwise} \\
    \end{cases}\notag\\
\rhos{4}(\P) &= \min_j \frac{\min_{\E\in\calI_\P^j}\mE}{\max_{\E\in\calI_\P^j}\mE}
\label{eq:indicator:2D}
\end{align}

For polyhedral meshes, we assess the interior quality of an element
$\P$ with a volumetric quality operator. 
Additionally, we
evaluate the quality of the element's faces, denoted as $\calF_{\P}$,
by employing the two-dimensional quality indicators $\rhos{i}$. 
This dual approach ensures a comprehensive analysis of
both the volumetric and surface attributes of the mesh.

\begin{align}
\rhov{1}(\P) &= \frac{\lvert kernel(\P)\rvert}{\mP}
    \prod_{\F\in\calF_{\P}}\rhos{1}(\F) \notag\\
    &= \begin{cases} 
    1 & \mbox{if $\P$ and all its faces are convex} \\
    \in (0,1) & \mbox{if $\P$ and all its faces are concave and star-shaped} \\
    0 & \mbox{if $\P$ or one of its faces are not star-shaped}
    \end{cases} \notag\\
\rhov{2}(\P) &= \frac{1}{2}\sqbra{
    \frac{\min(\sqrt[3]{\mP},\,\min_{\F\in\calF_{\P}}\hF)}{\hP} +
    \frac{\sum_{\F\in\calF_{\P}}\rhos{2}(\F)}{\#\cubra{\F\in\calF_{\P}}}} \notag\\
\rhov{3}(\P) &= \frac{1}{2}\sqbra{
    \frac{4}{\#\cubra{\F\in\calF_{\P}}} +
    \frac{\sum_{\F\in\calF_{\P}}\rhos{3}(\F)}{\#\cubra{\F\in\calF_{\P}}}}\notag\\
    &=\begin{cases} 
    1 & \mbox{if $\P$ is a tetrahedron} \\
    \in (0,1) & \mbox{otherwise} \\
    \end{cases}
    \label{eq:indicator:3D}
\end{align}

The $kernel$ operator in $\rhos{1}$ and $\rhov{1}$ computes the kernel
of an element (polygon or polyhedron), intended as the set of points
from which the whole element is visible.
The sets $\calI_\P^j$ in $\rhos{4}$ are the 1-dimensional disjoint
sub-meshes corresponding to the edges of $\P$ (we consider each
$\calI_\P^j$ as a mesh, as it may contain more than one edge) such
that $\calI_\P=\cup_j\calI_\P^j$, where $\calI_\P$ is the
1-dimensional mesh induced by $\partial\P$ (see
\cite{sorgente2021role}).
Indicator $\rhos{4}$ does not have a natural 3D extension~\cite{sorgente2021indicator}.

These indicators are combined together into an \textit{element quality
  indicator}, which measures the overall quality of an element
$\P\in\Th$:
\begin{align}
\varrho_2(\P) &= \sqrt{\frac{\rhos{1}(\P)\rhos{2}(\P) + \rhos{1}(\P)\rhos{3}(\P)+ \rhos{1}(\P)\rhos{4}(\P)}{3}}\notag\\
\varrho_3(\P) &= \sqrt{\frac{\rhov{1}(\P)\rhov{2}(\P) + \rhov{1}(\P)\rhov{3}(\P)}{2}}
\label{eq:indicator:element}
\end{align}
For $d=2,3$ we have $\varrho_d(\P)\approx 1$ if $\P$ is a
perfectly-shaped element, e.g. an equilateral triangle/tetrahedron or
a square/cube, $\varrho_d(\P)=0$ if and only if $\P$ is not
star-shaped, and $0<\varrho_d(\P)<1$ otherwise.

It is also possible to define a global function, which we call
\textit{mesh quality indicator} and denote by $\varrho_d (\Th)$ with a small
abuse of notation, to measure the overall quality of a $d$-dimensional mesh $\Th$:
\begin{equation}
\varrho_d(\Th) = \sqrt{\frac{1}{\#\left\{ \P \in \Th \right\}} \ \sum_{\P\in\Th} \varrho_d(\P)}.
\label{eq:indicator:mesh}
\end{equation}
All indicators in \eqref{eq:indicator:2D}, \eqref{eq:indicator:3D},
and consequently the local and global indicators
\eqref{eq:indicator:element}, \eqref{eq:indicator:mesh}, only depend
on the geometrical properties of the mesh elements; therefore their
values can be computed before applying the VEM, or any other numerical
scheme.

\subsection{Mesh Partitioning}
\label{subsec:optimization:metis}
METIS~\cite{karypis1997metis} is a collection of serial programs
implementing graph partitioning algorithms and finite element meshes,
and producing fill-reducing orderings for sparse matrices.
In the context of our work, we use it to perform a $K$-way
partitioning of the dual graph of a mesh.
We construct this graph by representing each mesh element as a node
and establishing arcs between adjacent elements, as depicted in
Figure~\ref{fig:pipeline:1}.
Additionally, we assign weights to both nodes and arcs based on
the element quality indicator \eqref{eq:indicator:element},
cf. Figure~\ref{fig:pipeline:2}.
In detail:
\begin{itemize}
\item Each mesh element $\P \in \Th$ becomes a graph
  node $n_{\P}$ with the associated weight $w(n_{\P})=\varrho_d(\P)$;
\item Each internal mesh $(d-1)$-dimensional object (which correspond
  to a face $\F \in \calF_{\hh}$ if $d = 3$ or an edge $e \in
  \calE_{\hh}$ if $d = 2$), shared by elements $\P_1$ and $\P_2$,
  becomes a graph arc $a_{\F}$ with weight
  $w(a_{\F})=\varrho_d(\P_1\cup\P_2)$;
\item $(d-1)$-dimensional boundary objects (boundary faces if $d=3$ or boundary edges if $d=2$) are not
  part of the graph, as they are only shared by one element.
\end{itemize}
We store the weights in two arrays $\boldsymbol{w_n}$, $\boldsymbol{w_a}$ that will be passed to METIS together with the mesh.
To maximize their impact on the mesh partition process, and also because METIS only accepts integer weights, we re-scale the weights in $\boldsymbol{w_n}$ from $[0,1]$ onto the interval 
\[
\left[ 1+\left\lfloor\min(\boldsymbol{w_n}) \frac{\#\boldsymbol{w_n}}{c_w}\right\rfloor, 
\ 1+\left\lfloor\max(\boldsymbol{w_n}) \frac{\#\boldsymbol{w_n}}{c_w}\right\rfloor \right]
\]
and similarly for $\boldsymbol{w_a}$. 
We add a constant 1 because zero-weights are not accepted, and divide by $c_w\in\REAL$ to avoid extremely big weights in very large meshes which may cause overflow problems.
For the meshes in this work we used $c_w=10$.

The routine \textit{METIS-PartGraphKway} partitions a graph into $K$ parts using multilevel $K$-way partitioning.
We compute the number of required partitions $K$ as a percentage of the number of elements in the mesh $\#\Th$.
This percentage is called the \textit{optimization parameter} $\calK$ and it is defined by:
\begin{equation*}
    K=\param\%(\#\Th):=\param\,\frac{\#\Th}{100}.
\end{equation*}
If we set $\param=20$ we are asking the optimization to partition the mesh into $K=20\%(\#\Th)$ parts, i.e., to generate a mesh with 1/5 of the elements of $\Th$.
In our experiments, we could observe that values of $\param$ below 5 lead to collapsing all the elements into a single one, while for $\param$ values higher than 40 the number of partitions computed by METIS does not increase proportionally.
These observations could highlight possible limitations of the partitioning algorithm that are independent of the quality indicator used.
The \textit{METIS-PartGraphKway} routine can be configured with a number of flags. In particular, we set \textit{PartitionType: CutBalancing} and activated the flags \textit{ContigousPartitions}, \textit{CompressGraph}, and \textit{MinimizeConnectivity} (see the METIS documentation for the details \cite{karypis1997metis}).

\subsection{Mesh Agglomeration}
\label{subsec:optimization:agglomeration}
The partition computed by METIS consists of an array containing a flag $\lambda_{\P}\in[1,K]$ for each element $\P$ of the mesh, see Figure~\ref{fig:pipeline:3}.
We aim to pick all the elements with the same flag and replace them with their union, which will be treated as a single element.
In particular, for each flag $\lambda=1,\ldots,K$ we consider the set $\Lambda:=\{\P\in\Th: \lambda_{\P}=\lambda\}$ and build a new element $\P_{\Lambda}$ defined by all the faces (edges if $d=2$) that are shared by only one element in $\Lambda$.
While doing this operation, we check that the arising element is still a manifold.
If not, we stop the agglomeration and skip to the next label, maintaining the old elements.
This situation is extremely rare, but a single non-manifold element may compromise the whole mesh for certain applications.

At the end of the optimization process with input parameter $\param$, the resulting mesh will have approximately $K$ elements optimized for the computation of the virtual element basis functions.
We will have removed around $\#\Th-K$ elements from the mesh and at least $\#\Th-K$ faces (edges if $d=2$) between them, see Figure~\ref{fig:pipeline:4}.
Moreover, if $K$ is particularly low, we may happen to remove also edges (if $d=3$) and vertices, but in a smaller number.

\section{Datasets}
\label{sec:datasets}
This section details the generation and the optimization of a number of meshes specifically created for testing our algorithm.
More precisely, we call \textit{dataset} a sequence of meshes defined over the same domain, with decreasing mesh size, and built with the same technique, so that they contain similar elements organized in similar configurations.
We generate four datasets to span the most common types of meshes currently used in numerical simulations: triangular and quadrangular in 2D, and tetrahedral and hexahedral in 3D.
These meshes are built with random strategies so that they contain low-quality elements that make it meaningful, if not necessary, to think about their optimization.
All the meshes presented in this paper are available for download at \url{GitHub} \textit{(GitHub link to be added before publication)}.

\subsection{Planar Datasets}
\label{subsec:datasets:planar}
\begin{figure}[!h]
  \centering
  \begin{subfigure}[]{.3\textwidth}
      \centering
    \includegraphics[width=\linewidth]{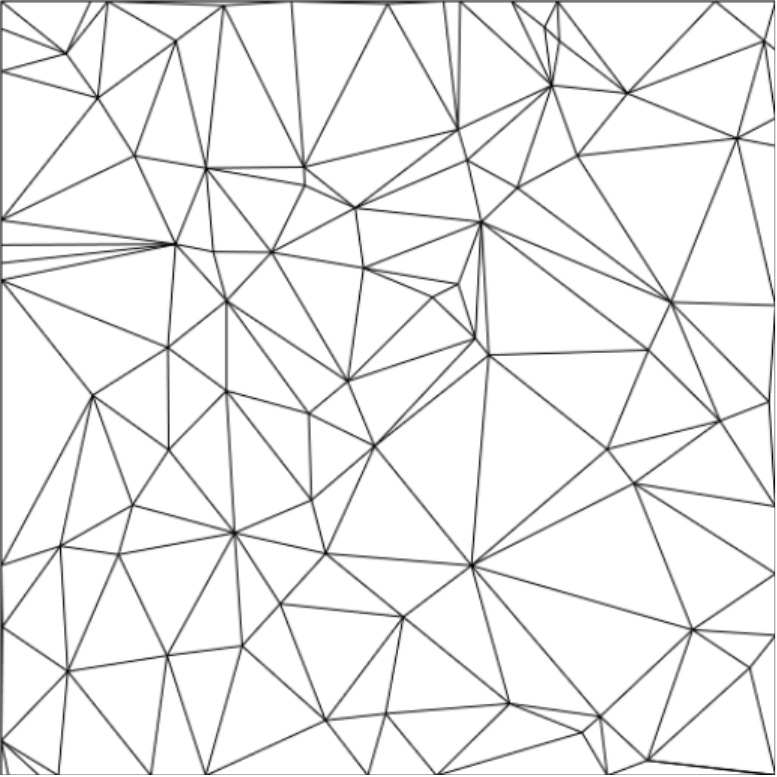}
      \caption{$\dtri{}$}
      \label{fig:dataset:tri100}
  \end{subfigure}
  \hspace{.2cm}
  \begin{subfigure}[]{.3\textwidth}
      \centering
    \includegraphics[width=\linewidth]{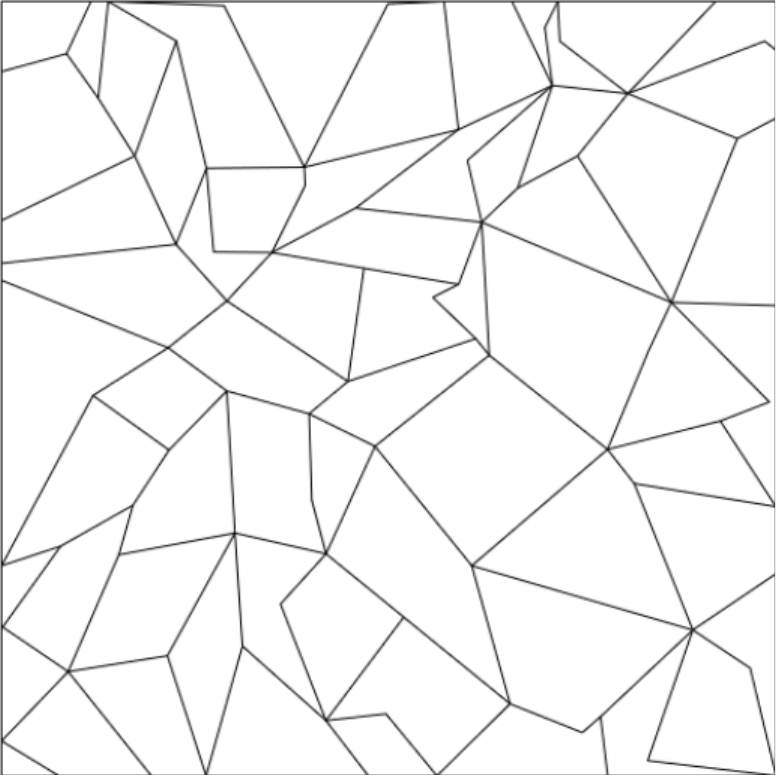}
      \caption{$\dtri{40}$}
      \label{fig:dataset:tri40}
  \end{subfigure}
  \hspace{.2cm}
  \begin{subfigure}[]{.3\textwidth}
      \centering
    \includegraphics[width=\linewidth]{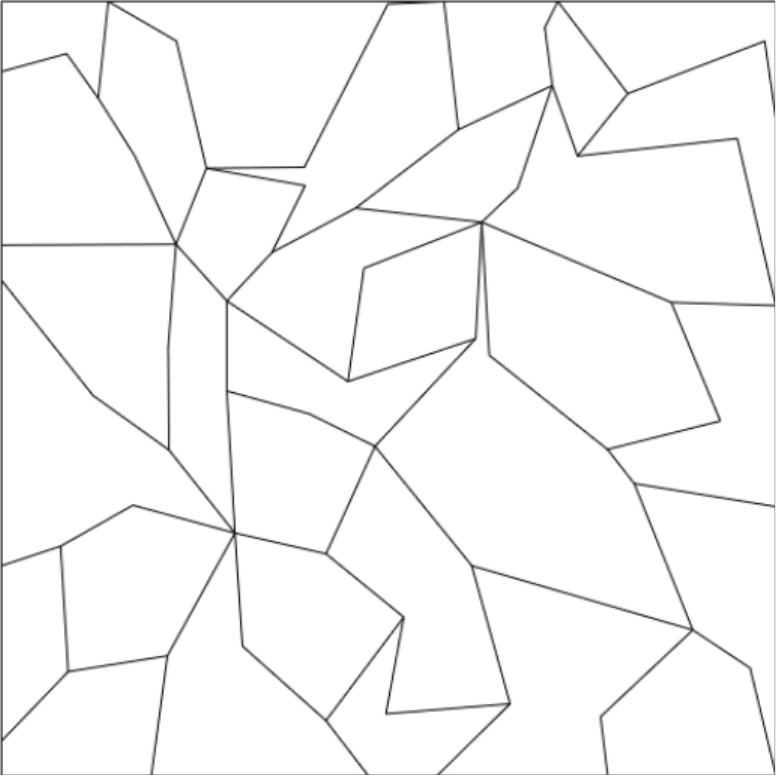}
      \caption{$\dtri{20}$}
      \label{fig:dataset:tri20}
  \end{subfigure}
  
  \vspace{.5cm}
  \begin{subfigure}[]{.3\textwidth}
      \centering
    \includegraphics[width=\linewidth]{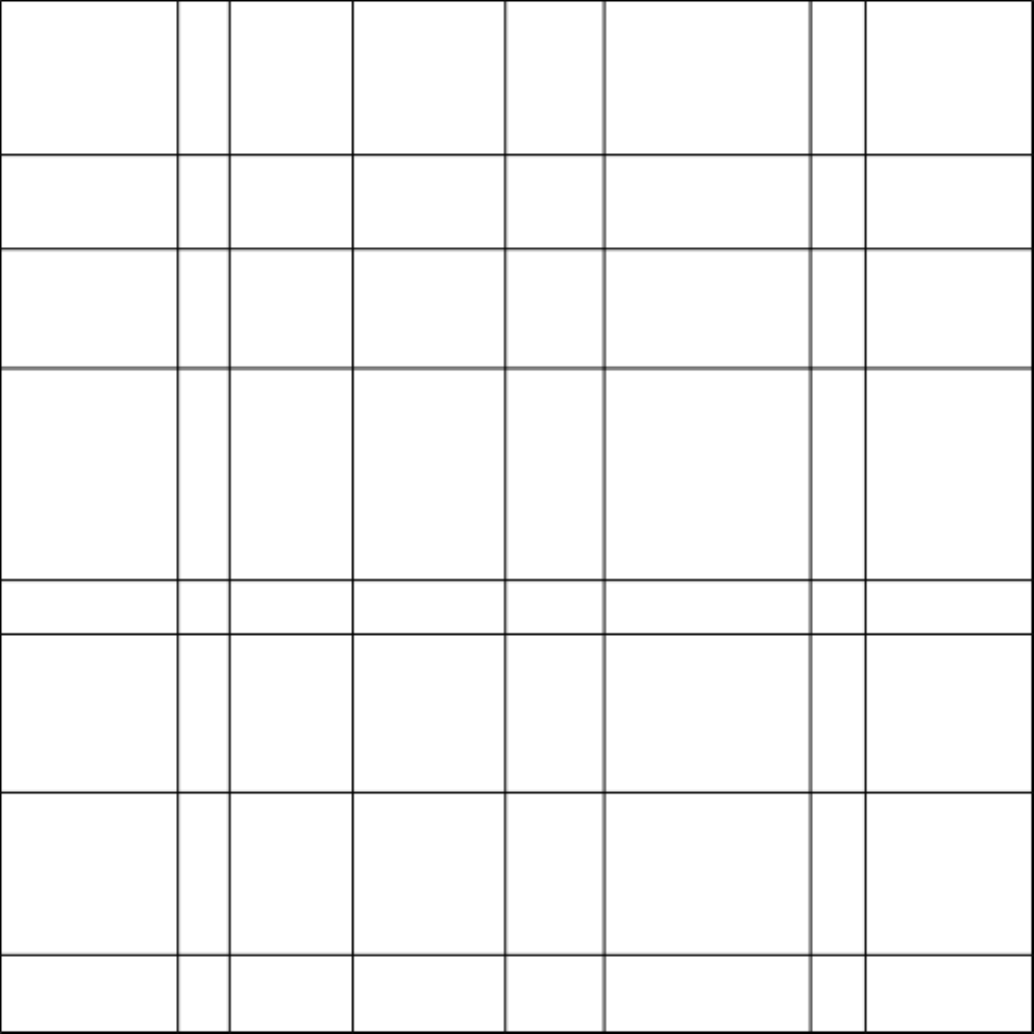}
      \caption{$\dqua{}$}
      \label{fig:dataset:quad100}
  \end{subfigure}
  \hspace{.2cm}
  \begin{subfigure}[]{.3\textwidth}
      \centering
    \includegraphics[width=\linewidth]{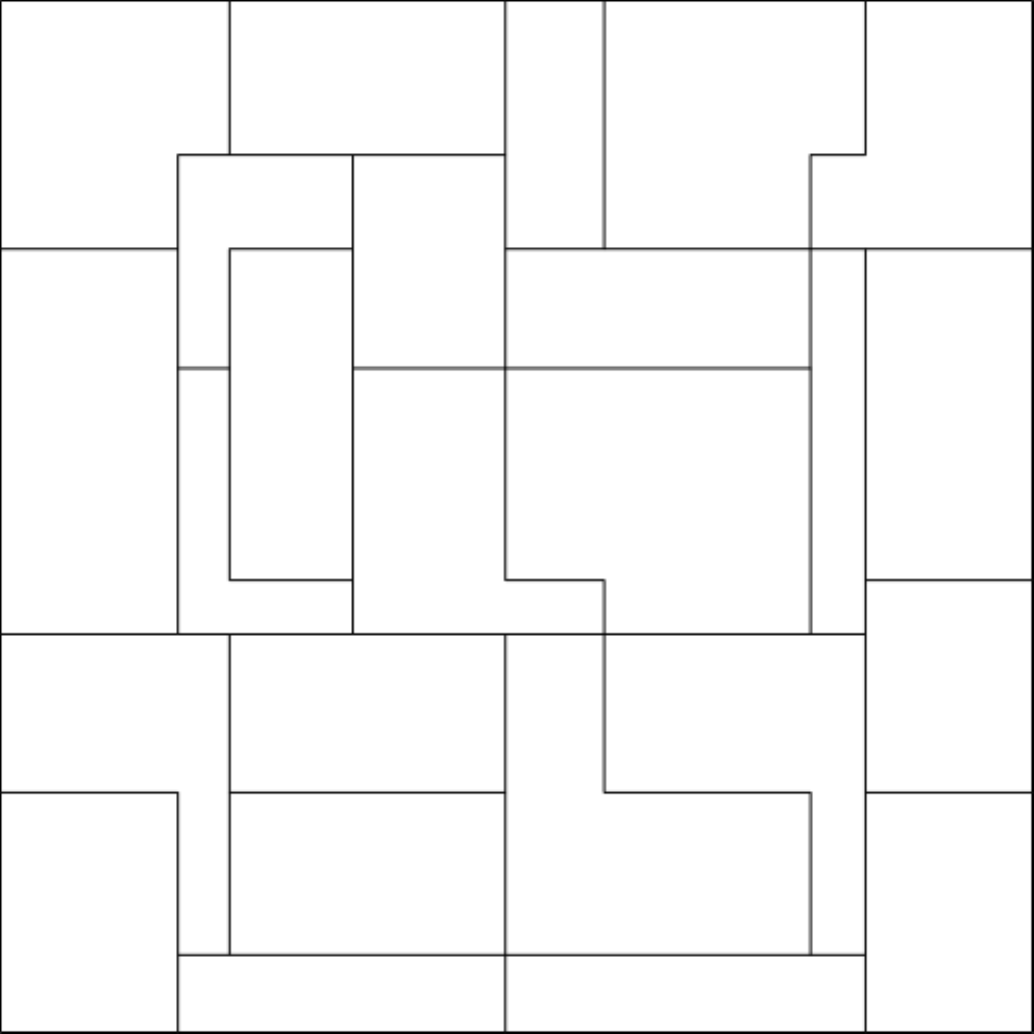}
      \caption{$\dqua{40}$}
      \label{fig:dataset:quad40}
  \end{subfigure}
  \hspace{.2cm}
  \begin{subfigure}[]{.3\textwidth}
      \centering
    \includegraphics[width=\linewidth]{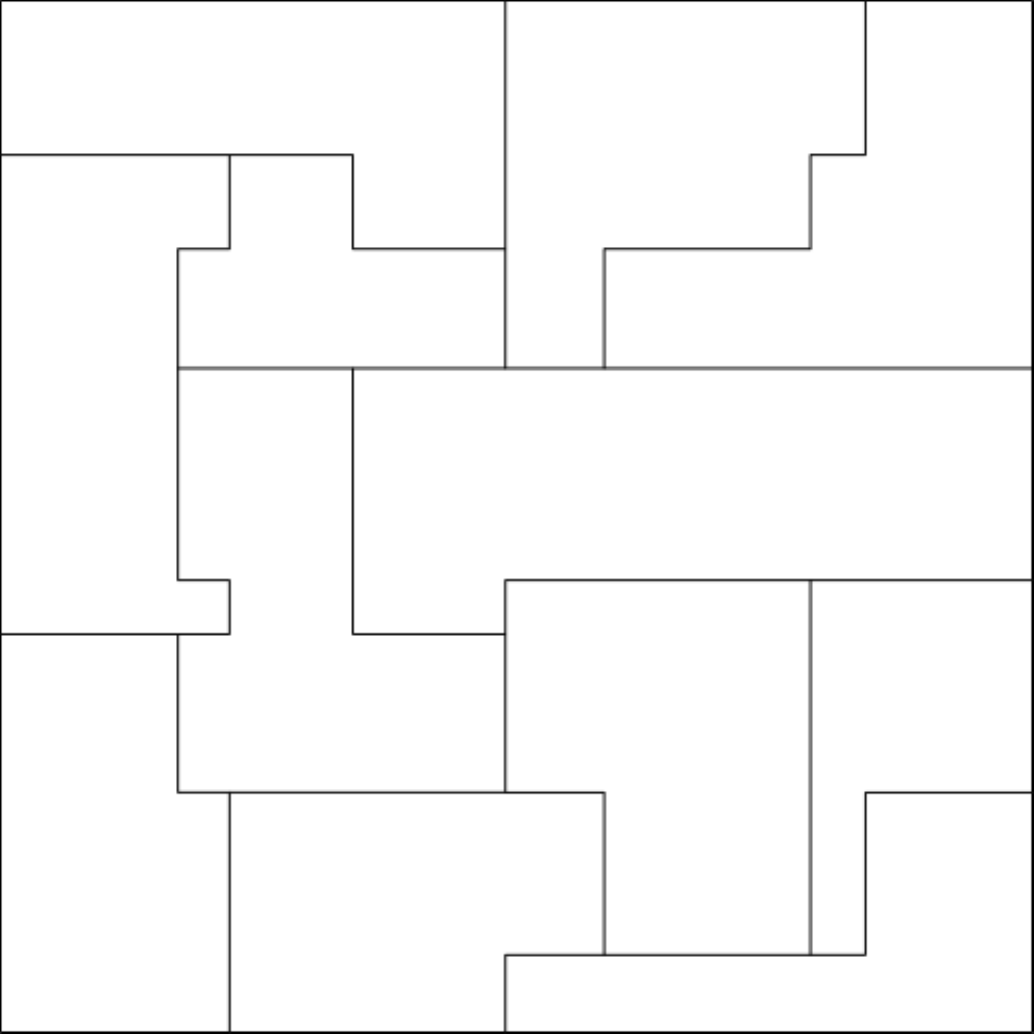}
      \caption{$\dqua{20}$}
      \label{fig:dataset:quad20}
  \end{subfigure}
  \caption{First mesh of datasets $\dtri{}$ and $\dqua{}$ optimized with $\param=40$ and $\param=20$.}
  \label{fig:dataset:2D}
\end{figure}

We generate two planar datasets on the unit square $\Omega=\left[0,1\right]^2$, containing five meshes each.
To compare our results, the $i-$th mesh of each dataset contains a similar number of vertices and elements.

Dataset $\dtri{}$ contains triangular meshes with randomly-shaped elements, generated by calling \textit{Triangle}~\cite{shewchuk1996triangle} on $\Omega$ with a set of initial points to be preserved and no other constraints.
The initial points are randomly sampled on $\Omega$. Hence, the resulting meshes contain several needles and flat elements.
Dataset $\dqua{}$ contains quadrangular meshes with randomly-shaped quads.
We start from an equispaced planar grid and randomly move all the points belonging to a certain plane $\pi_1$ to another plane $\pi_2$, parallel to $\pi_1$.
In practice, we pick all the points sharing the same $x$-coordinate and randomly change $x$ by the same quantity; then, we repeat this operation for the $y$ coordinates.
The resulting meshes contain quadrangular elements with very different sizes and aspect ratios.
In Figure~\ref{fig:dataset:2D} we present the first mesh of datasets $\dtri{}$ and $\dqua{}$, and their optimized versions with parameters $\param=40$ and $\param=20$.

In the datasets optimized with parameter $40$, each mesh has less than half the elements of its corresponding original mesh, while in datasets optimized with parameter $20$ each mesh has around $1/5$ of the elements of its corresponding original mesh.
Moreover, while in the original datasets all the elements are guaranteed to be convex (triangles or quads), in the optimized ones a significant number of elements are not convex or even not star-shaped.

\subsection{Volumetric Datasets}
\label{subsec:datasets:volumetric}
\begin{figure}[!h]
  \centering
  \begin{subfigure}[]{.3\textwidth}
      \centering
    \includegraphics[width=\linewidth]{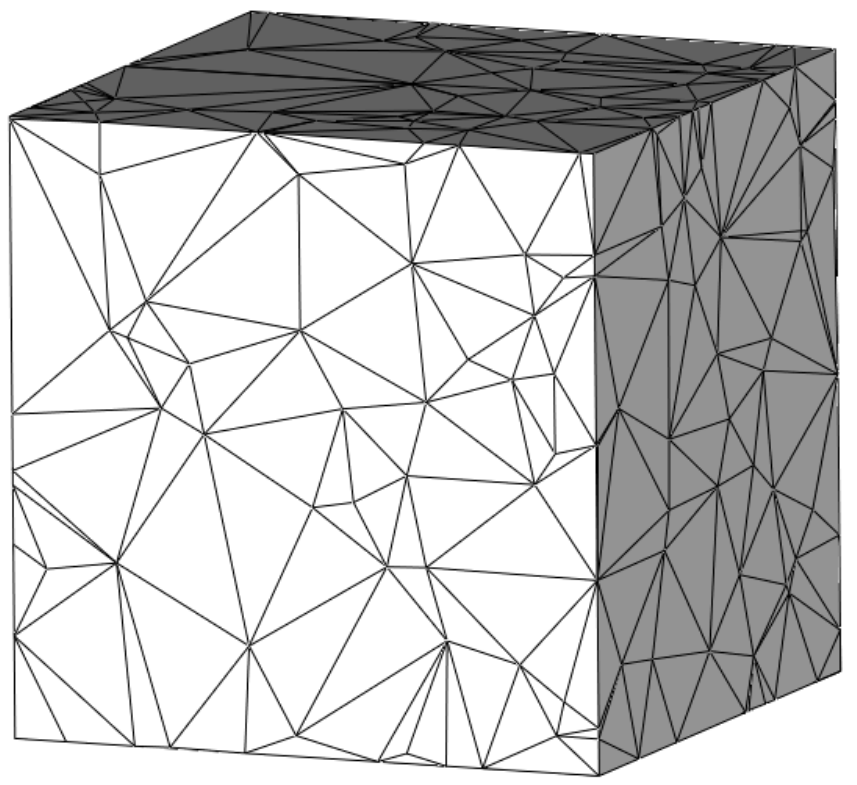}
      \caption{$\dtet{}$}
      \label{fig:dataset:tet100}
  \end{subfigure}
  \hspace{.2cm}
  \begin{subfigure}[]{.3\textwidth}
      \centering
    \includegraphics[width=\linewidth]{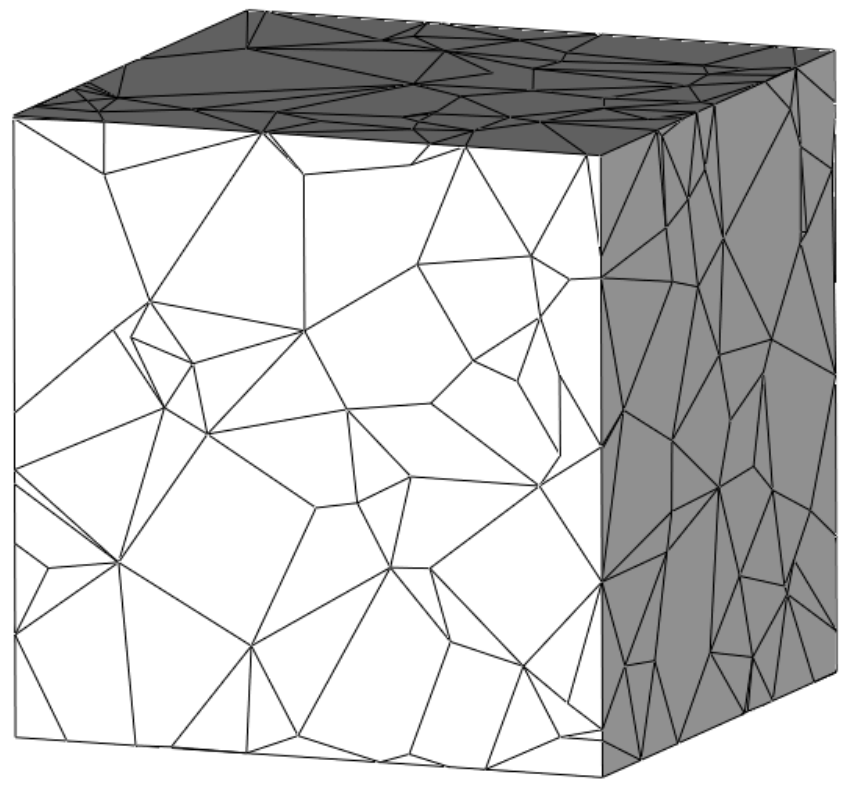}
      \caption{$\dtet{40}$}
      \label{fig:dataset:tet40}
  \end{subfigure}
  \hspace{.2cm}
  \begin{subfigure}[]{.3\textwidth}
      \centering
    \includegraphics[width=\linewidth]{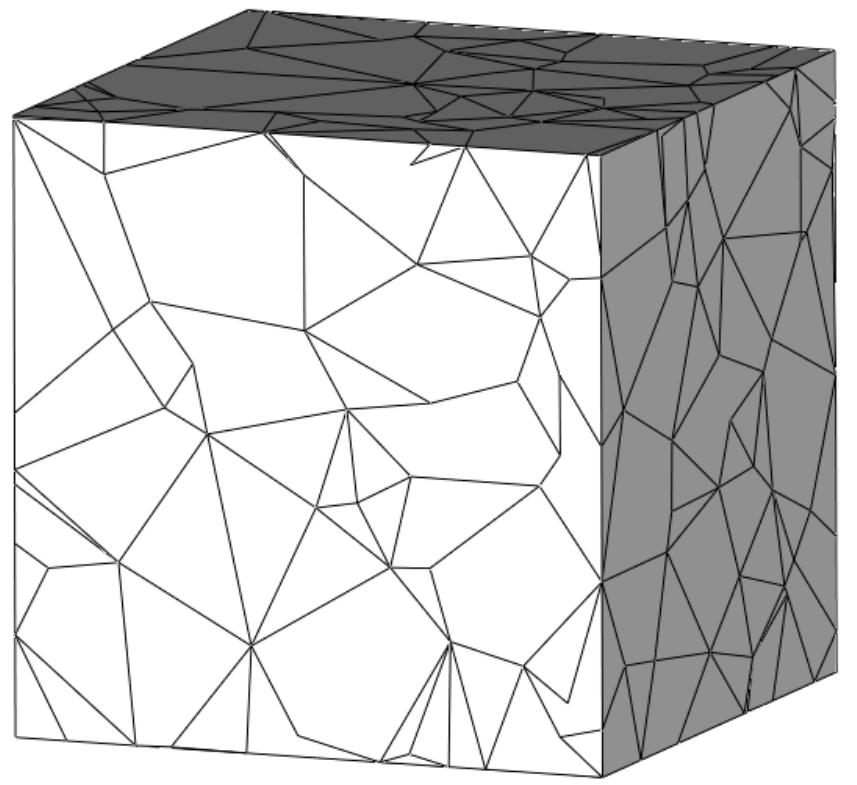}
      \caption{$\dtet{20}$}
      \label{fig:dataset:tet20}
  \end{subfigure}
  
  \vspace{.5cm}
  \begin{subfigure}[]{.3\textwidth}
      \centering
    \includegraphics[width=\linewidth]{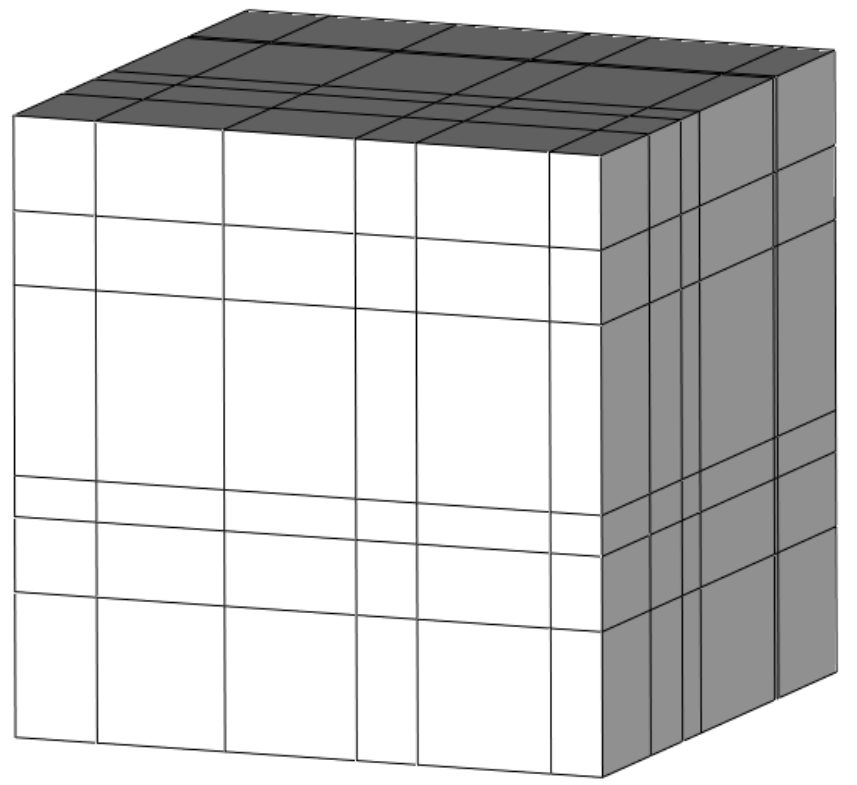}
      \caption{$\dhex{}$}
      \label{fig:dataset:hex100}
  \end{subfigure}
  \hspace{.2cm}
  \begin{subfigure}[]{.3\textwidth}
      \centering
    \includegraphics[width=\linewidth]{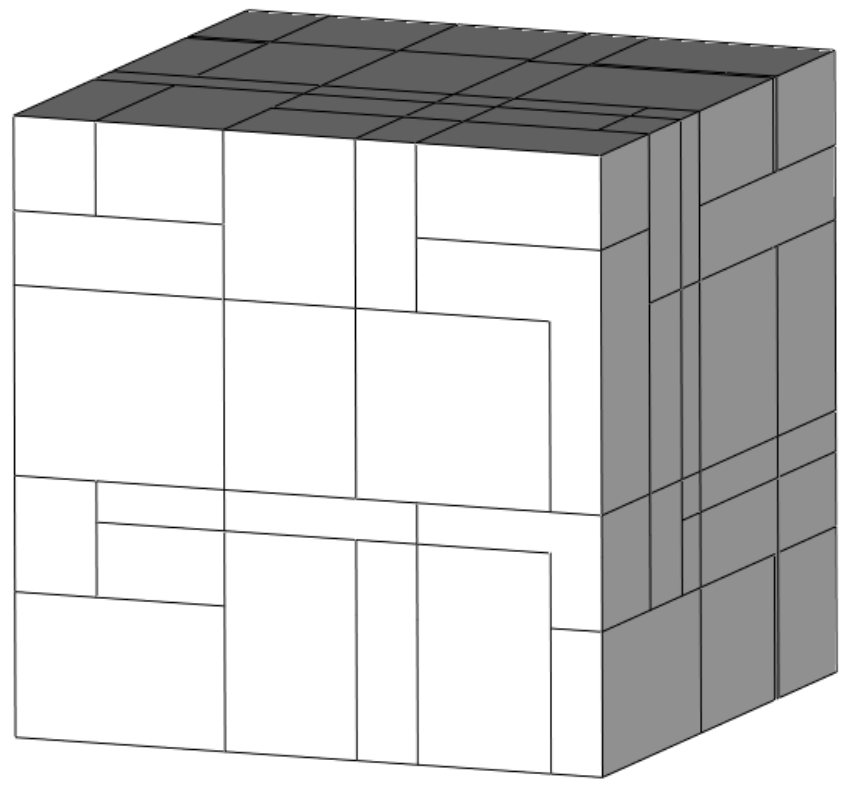}
      \caption{$\dhex{40}$}
      \label{fig:dataset:hex40}
  \end{subfigure}
  \hspace{.2cm}
  \begin{subfigure}[]{.3\textwidth}
      \centering
    \includegraphics[width=\linewidth]{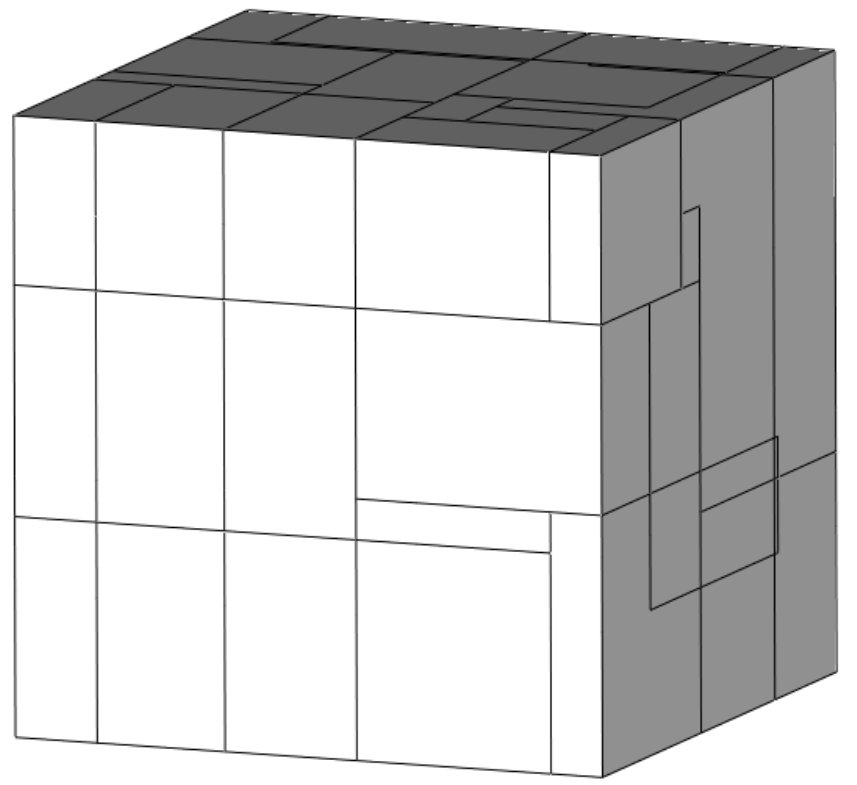}
      \caption{$\dhex{20}$}
      \label{fig:dataset:hex20}
  \end{subfigure}
  \caption{First mesh of datasets $\dtet{}$ and $\dhex{}$ optimized with $\param=40$ and $\param=20$.}
  \label{fig:dataset:3D}
\end{figure}
\begin{table}[!h]
\centering
\caption{Geometric properties of datasets $\dtet{}$, $\dtet{40}$, and $\dtet{20}$:
number of vertices, faces and elements of each mesh; number of internal DOFs for $k=1,2,3$; mesh quality~\eqref{eq:indicator:mesh}.} \label{tab:optimization}
\includegraphics[width=.9\textwidth]{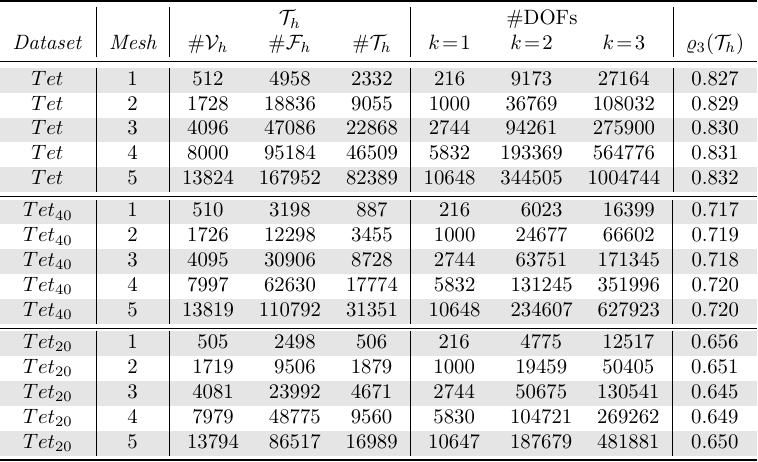} 
\end{table}

We also define two volumetric datasets on the unit cube $\Omega=\left[0,1\right]^3$, containing five meshes each, with similar numbers of vertices and elements.
They are built in analogy to their planar versions: for $\dtet{}$ we use \textit{TetGen}~\cite{si2006quality} and for $\dhex{}$ we also perturb the $z$-coordinates.
In Figure~\ref{fig:dataset:3D} we show the first mesh of these datasets and their optimized versions with parameters $\param=40$ and $\param=20$.

In Table~\ref{tab:optimization}, we present some numerics about dataset $\dtet{}$ and its optimized versions; analogous tables for the other datasets are reported in~\ref{app} (Tables~\ref{app:dataset:optimization:2D}, \ref{app:dataset:optimization:3D}).
As already noted in Section~\ref{subsec:optimization:agglomeration}, the number of vertices $\#\calV_h$ does not change significantly after the optimization process.
Instead, the number of elements $\#\Th$ decreases as expected: meshes in $\dtet{40}$ and $\dtet{20}$ contain, approximately, $40\%$ and $20\%$ the elements of the corresponding meshes in $\dtet{}$.
The number of faces $\#\calF_h$ decreases a bit more than the number of elements (in absolute value), given that each couple of merged elements shares at least one face which gets erased.
For instance in mesh 5 from $\dtet{}$ to $\dtet{40}$ we remove 51038 elements and 57160 faces, obtaining a mesh with $38\%$ of the elements and $66\%$ of the faces of the original one.

These considerations are reflected in the last three columns of the tables, in which we report the number of internal degrees of freedom, for $k=1,2,3$, that we would define over each mesh in a VEM simulation.
For $k=1$ the DOFs coincide with the internal vertices and therefore remain almost identical, but for $k>1$ we have an important reduction.
For $\dtet{40}$ we save on average the $33\%$ in DOFs for $k=2$ and the $38\%$ for $k=3$, while for $\dtet{20}$ we save on average the $47\%$ in DOFs for $k=2$ and the $53\%$ for $k=3$.
This means that, for instance, optimizing mesh 5 with $\param=20$ we can save up to $5\times10^6$ DOFs for $k=3$, obtaining a mesh that is significantly lighter to store and easy to handle.

We point out that the primary goal of the algorithm is to reduce the number of elements to the expected percentage.
In this sense, the global mesh quality is likely to decrease after the optimization, as visible in the last column of Table~\ref{tab:optimization}.
However, the loss of quality is optimal with respect to the number of elements removed.
In other words, the optimized mesh has the maximum quality we can achieve by removing that number of elements.

\section{Convergence Tests}
\label{sec:convergence}
In this section, we compare the performance of the VEM over the original datasets and the optimized ones.

\subsection{Test Problem}
\label{subsec:convergence:problem}
On all the datasets from Section~\ref{sec:datasets} we solve the Poisson problem:
\begin{align}
    - \Delta u &= f\phantom{0} \text{in } \Omega,\label{eq:prob:1:A} \\
    u          &= 0\phantom{f} \text{in } \partial \Omega,\label{eq:prob:1:B}\nonumber
\end{align} 
reported in the variational form of Problem~\eqref{eq:prob:var} with the definition of
\begin{equation*}
    a(u, v) := \left( \nabla u, \nabla v \right)_{\Omega},\ F(v) := \left( f, v \right)_{\Omega}.
\end{equation*}
For the discretization, we employ the VEM described in Section~\ref{sec:problem}.
In particular, we define the approximated bilinear form of Equation~\eqref{eq:vem:a} as:
\begin{equation*}
    a^{\P}_{h}(u_h, v) 
    = (\Pi^0_{k - 1} \nabla u_h, \Pi^0_{k - 1} \nabla v)_{\Omega}
    + S^{\P}(u - \Pi^{\nabla}_k u, v - \Pi^{\nabla}_k v).
\end{equation*}
and $F_h^{\P}(v) := (f, \Pi^0_k v)_{\Omega}$.
We use the symbol $\Pi^0_{k-1}$ for the $L^2$-projection operator of a vector-valued function onto the polynomial space, applied component-wisely.
Finally, we compute $f \in L^2(\Omega)$, using as ground truth the function 
\begin{equation*}
    u(x, y, z) = 6.4 xyz(x-1)(y-1)(z-1).
\end{equation*}
We compare the results plotting the relative $\LTWO$-norm and $\HONE$-seminorm:
\begin{align}
    \errL &= \norm{u-u_h}{0,\Omega}\,/\,\norm{u}{0,\Omega}, \label{eq:err_L2} \\
    \errH &= \snorm{u-u_h}{1,\Omega}\,/\,\snorm{u}{1,\Omega}, \label{eq:err_H1} 
\end{align}
of the approximation error $u-u_h$ as the number of degrees of freedom increases.
The optimal convergence rate of the method is indicated by the slope of a reference triangle in each plot.
We want to ensure that the method converges with the correct rate over the optimized datasets, and we are interested in measuring differences in the approximation errors produced before and after the optimization.

\medskip
Given a mesh $\Th$ over which we have computed an approximated solution $u_h$ of $u$, and an optimized mesh $\Th'$ with approximated solution $u_h'$, we measure the difference between $\Th$ and $\Th'$ through the following quantities:
\begin{align} 
    \delD &= \#\text{dofs}(\Th) - \#\text{dofs}(\Th'), \label{eq:delta_dofs} \\
    \delL &= \frac{\norm{u-u_h}{0,\Omega} - \norm{u-u_h'}{0,\Omega}}{\norm{u}{0,\Omega}}, \label{eq:delta_L2} \\
    \delH &= \frac{\snorm{u-u_h}{1,\Omega} - \snorm{u-u_h'}{1,\Omega}}{\snorm{u}{1,\Omega}}. \label{eq:delta_H1}
\end{align}
In practice, $\delL>0$ means that the relative error produced by the VEM over $\Th'$ is smaller than the one produced over $\Th$, and similarly for $\delH$ and $\delD$.
We expect optimized meshes to produce less accurate results than the original ones (i.e., $\delL<0$ and $\delH<0$) because these discretizations contain many fewer degrees of freedom and, therefore, information.
In fact, $\delD\ge0$ for all meshes.
However, if the loss in accuracy is limited, the optimization might be convenient.
For instance, if $\delL\sim~-10^{-7}$ and $\delH\sim~-10^{-4}$, but $\delD\sim~10^{5}$, we have a significant gain in terms of space and computations at a very small cost in terms of accuracy.

Finally, we compare the computational time:
\begin{align} 
    \delT &= T_{\text{solve}}(\Th) - (T_{\text{optimize}}(\Th) + T_{\text{solve}}(\Th')), \label{eq:delta_t}
\end{align}
where $T_{\text{solve}}(\Th)$ is the time required for solving \eqref{eq:prob:1:A} on $\Th$, and $T_{\text{optimize}}(\Th)$ is the time required for optimizing $\Th$ into $\Th'$.

\subsection{Planar Datasets}
\label{subsec:convergence:planar}
\begin{figure}[!h]
  \centering
  \includegraphics[width=\textwidth]{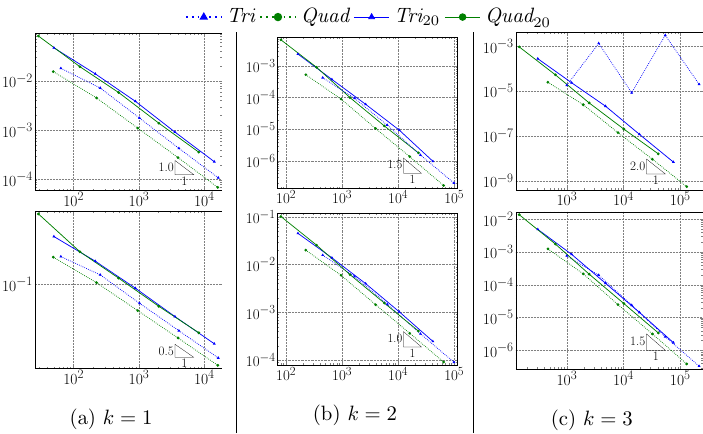}
  \caption{Convergence of datasets $\dtri{20}$, $\dqua{20}$ (continuous lines) compared to their original datasets $\dtri{}$, $\dqua{}$ (dotted line). We measure $\errL$ (top figures) and $\errH$ (bottom figures) with respect to the number of DOFs.}
  \label{fig:results:planar}
\end{figure}

We begin our analysis by checking the convergence of the VEM for $k=1,2,3$ over the planar datasets $\dtri{}$ and $\dqua{}$ from Section~\ref{subsec:datasets:planar} and their optimizations with $\param=20$.
First, we note from the convergence plots in Figure~\ref{fig:results:planar} how the optimized datasets $\dtri{20}$ and $\dqua{20}$ preserve the optimal convergence rate both in $\LTWO$ and $\HONE$ norms.

We then consider the quantities~\eqref{eq:delta_dofs}, \eqref{eq:delta_L2}, and \eqref{eq:delta_H1}, computed between the last meshes (mesh $5$) of $\dtri{}$ and $\dtri{20}$:
\begin{itemize}
    \item for $k=1$, we have $\delD\sim~10^3$, $\delL\sim~-10^{-4}$, and $\delH\sim~-10^{-3}$;
    \item for $k=2$, we have $\delD\sim~10^4$, $\delL\sim~-10^{-7}$, and $\delH\sim~-10^{-4}$;
    \item for $k=3$, we have $\delD\sim~10^5$, $\delL\sim~+10^{-5}$, and $\delH\sim~-10^{-6}$.
\end{itemize}
Full details for the other meshes are reported in Table~\ref{app:results:2D} of the Appendix.
We note how for $k=1,2$ the method performs slightly worse on the optimized meshes, despite having significantly fewer degrees of freedom, and this difference decreases as $k$ grows.
The difference in $\delD$ is appreciable because the dots relative to optimized meshes are shifted towards the left in Figure~\ref{fig:results:planar}.
For $k=3$ the error $\errL$ produced by the original meshes is so high, not only in the last mesh, that the method fails to converge.
This is due to extremely small and flat triangles in the original meshes, which cause enormous condition numbers.
Such elements disappear in the optimized meshes, and the VEM converges properly over them.
For instance, in dataset $\dtri{}$ we have $cond\sim~10^{12}$ and $cond\sim~10^{15}$ for meshes $3$ and $4$, while in $\dtri{20}$ we have $cond\sim~10^{5}$ and $cond\sim~10^{6}$.
Therefore, we have fewer DOFs and smaller $\errL$ in this case.

For $\dqua{}$ and $\dqua{20}$, the values of $\delD$, $\delL$, and $\delH$ are essentially analogous, see Table~\ref{app:results:2D}. 
The only significant difference is that the method does not diverge over $\dqua{}$ for $k=3$: in mesh 5 we have $\delD\sim~10^4$, $\delL\sim~-10^{-8}$, and $\delH\sim~-10^{-6}$.
The difference $\delL$ is negative in this case, but it is extremely small compared to $\delD$.
The difference in computational time between the original and the optimized meshes is negligible both for $\dtri{}$ and $\dqua{}$, as the considered meshes are relatively small.
The quantity $\delT$ will become significant with the more complex meshes of the next section.

\subsection{Volumetric Datasets}
\label{subsec:convergence:volumetric}
\begin{figure}[!h]
  \centering
  \includegraphics[width=\textwidth]{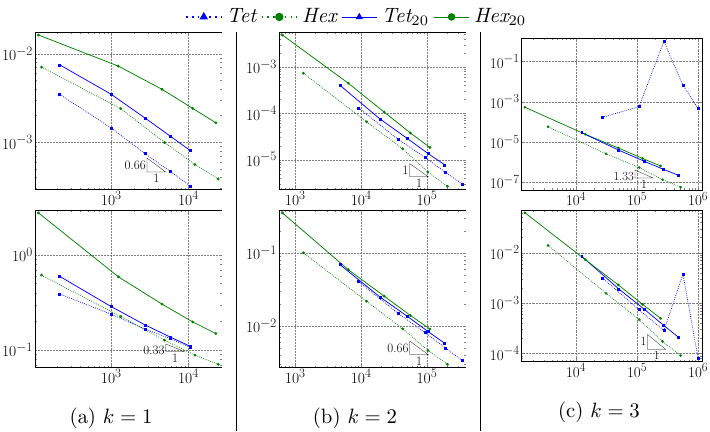}     
  \caption{Convergence of datasets $\dtet{20}$ and $\dhex{20}$ (continuous lines) compared to their original datasets (dotted lines). We measure $\errL$ (top figures) and $\errH$ (bottom figures) with respect to the number of DOFs.}
  \label{fig:results:volumetric}
\end{figure}
\begin{table}[!h]
  \centering
  \caption{VEM performance: difference between the number of DOFs~\eqref{eq:delta_dofs}, the $\errL$ error~\eqref{eq:delta_L2}, the $\errH$ error~\eqref{eq:delta_H1}, and the computational time~\eqref{eq:delta_t} w.r.t. the original mesh.} \label{tab:results}
  \includegraphics[width=.8\textwidth]{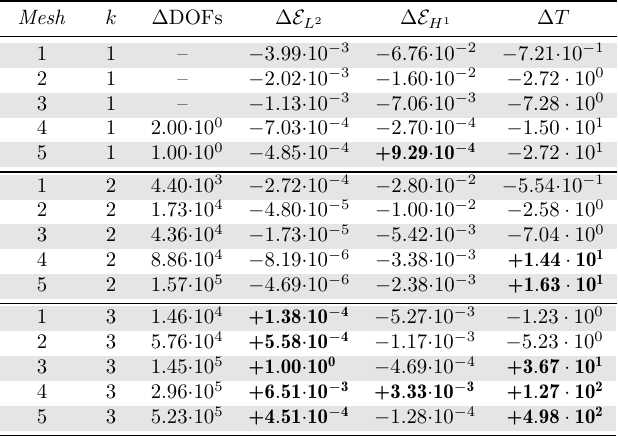}
\end{table}

In Figure~\ref{fig:results:volumetric} we present analogous plots for the 3D datasets $\dtet{}$ and $\dhex{}$, both optimized with $\param=20$.
As for the planar datasets, all optimized datasets produce optimal convergence rates, even when the original dataset does not (see $\dtet{}$ for $k=3$).

In Table~\ref{tab:results} we report the usual data about the differences between $\dtet{}$ and $\dtet{20}$:
\begin{itemize}
    \item for $k=1$ we have almost no DOFs reduction, with $\delL\sim~-10^{-4}$ and $\delH$ increasing until becoming a positive $10^{-4}$ in the last mesh;
    \item for $k=2$ $\delD$ reaches $10^5$ with $\delL\sim~-10^{-6}$ and $\delH\sim~-10^{-3}$;
    \item for $k=3$ the VEM diverges on the original meshes and we get $\delL>0$ on every mesh.
\end{itemize}
Results for $\dhex{}$ and $\dhex{20}$ can be found in Table~\ref{app:results:3D}.
For $k=1$ the difference between the original and the optimized meshes is higher, also in terms of $\delD$, but for $k>1$ the errors become smaller and smaller, even if the original dataset does not diverge in this case.

Looking at the times, for $k=1$ we have a negative $\delT$ for both $\dtet{20}$ and $\dhex{20}$.
This is because the elements of the optimized meshes generally have more complicated shapes than the simple tets or hexes of the original ones.
At the same time, for $k=1$ the number of DOFs does not decrease significantly, hence the linear system to be solved is more complex and has the same size.
However, once $k$ increases and $\delD$ becomes significant, the size of the system reduces drastically, and we see positive $\delT$ values for $k=2,3$.

In conclusion, we note how the optimization becomes more interesting for high values of $k$ and more complicated meshes.
An optimized mesh containing $10^5$ less DOFs, which produces results only $10^{-6}$ less accurate, leads to a faster and cheaper simulation.
Moreover, the optimization removes small and flat elements in meshes from $\dtet{}$ that cause the VEM to diverge, and we obtain a numerical approximation that was impossible to compute on the original meshes.

\subsection{Role of the Optimization Parameter}
\label{subsec:convergence:parameter}
\begin{table}[!h]
  \centering
  \caption{Role of parameter $\param$ in the optimization of dataset $\dtet{}$. 
  Each column shows the average (among the five meshes in the dataset) percentage reduction with respect to the original dataset.}
  \label{tab:results:parameter}
  \includegraphics[width=.48\textwidth]{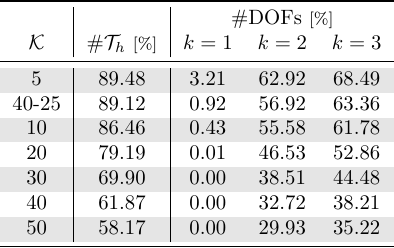}
\end{table}
\begin{figure}[!h]
  \centering
  \includegraphics[width=\textwidth]{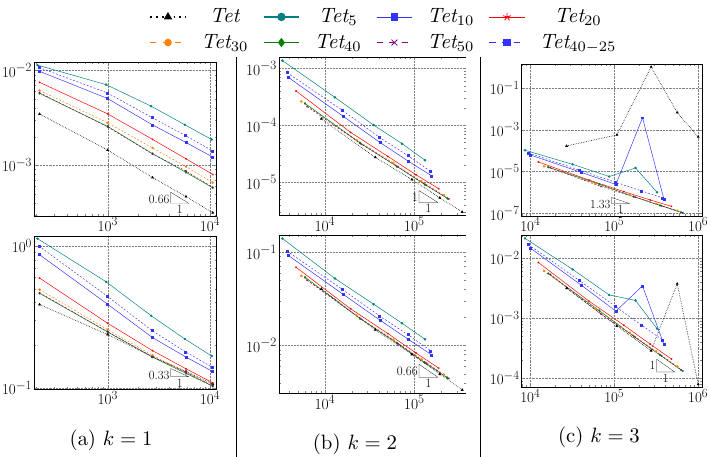}
  \caption{Comparison between dataset $\dtet{}$ and its optimized versions with $\param \in \{5, 10, 20, 30, 40, 50, \text{40-25}\}$. We measure $\errL$ (top figures) and $\errH$ (bottom figures) concerning the number of DOFs.}
  \label{fig:results:parameter}
\end{figure}

We now analyze the impact of parameter $\param$ on the optimization process, comparing the performance produced by dataset $\dtet{}$ optimized with $\param=5, 10, 20, 30, 40$, and 50.
We point out that $\dtet{}$ and $\dtet{20}$ are the same datasets shown in Figure~\ref{fig:results:volumetric}, and that we reported in Table~\ref{tab:optimization} the number of vertices, faces, elements and DOFs for meshes in $\dtet{}$, $\dtet{40}$, and $\dtet{20}$.
We also try to process the dataset recursively, i.e., to further optimize $\dtet{40}$ with $\param=25$, obtaining a dataset $\dtet{40-25}$.
This dataset contains the $25\%$ of the elements of $\dtet{40}$, which contains the $40\%$ of the elements of $\dtet{}$.
Overall, $\dtet{40-25}$ should be comparable with $\dtet{10}$.

In Table~\ref{tab:results:parameter} we present the reduction, in terms of numbers of elements $\#\Th$ and DOFs, of each optimized dataset compared to $\dtet{}$, and in Figure~\ref{fig:results:parameter} the VEM performance over them.
First, we note how the algorithm is generally able to satisfy the elements reduction required by the parameter: for $\param=40$ we have around $-60\%$ elements, for $\param=30$ we have around $-70\%$ elements, and so on.
As already noted in Section~\ref{subsec:optimization:agglomeration}, we reach a plateau for $\param>40$ because METIS implementation, therefore results for $\param=50$ are similar to those for $\param=40$.
At the opposite side of the range, when optimizing with $\param\le10$ the correspondence between $\param$ and the elements reduction becomes weaker ($86\%$ elements reduction instead of $90\%$).
However, if we subdivide the optimization into two steps, i.e. with $\param=40-25$, we get closer to the desired reduction.
Concerning the DOFs, for $k=1$ we have only a small reduction for the smallest $\param$ values, as expected.
For $k=2,3$ the DOFs reduction is smaller than the elements reduction, but they vary similarly.

In the plots of Figure~\ref{fig:results:parameter} we can observe how higher $\param$ values produce higher errors, but also that this difference decreases when increasing the order $k$.
Datasets $\dtet{10}$ and $\dtet{40-25}$ perform very similarly, but for $k=3$, the former starts to oscillate while the latter does not.
Hence, a single, aggressive optimization produces worse results than two conservative ones.

In conclusion, the choice of the best optimization parameter may depend on the order of the method.
In any case, we suggest to set $\param<50$ because the optimization does not work properly after that value.
For low-order VEM, it seems preferable to use a higher optimization value because removing too many elements may significantly increase the errors.
For $k>1$ instead, the difference between $\param=40$, $\param=30$, and $\param=20$ are very small in terms of $\errL$ while notable in the number of DOFs. Therefore, it is worth choosing a lower optimization value. 
Values of $\param$ lower than 20 are not recommended because they do not always produce reliable results, especially for $k>2$.
In those cases, applying multiple optimizations with higher $\param$ values is preferable.

\subsection{Importance of Quality Weights}
\label{subsec:convergence:weights}
We also analyze the importance of using the quality indicator~\eqref{eq:indicator:element} for assigning METIS weights to nodes and arcs (see Section~\ref{subsec:optimization:metis}).
In Figure~\ref{fig:results:weights} we consider dataset $\dhex{}$ optimized always with $\param=20$, but with different weights:

\begin{figure}[!h]
  \centering
  \includegraphics[width=\textwidth]{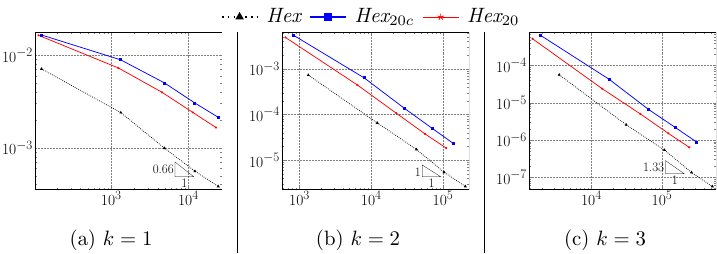}
  \caption{Error $\errL$ with respect to the number of DOFs produced by $\dhex{}$ optimized with $\param=20$ and constant weights ($\dhex{20c}$), or quality weights ($\dhex{20}$).}
  \label{fig:results:weights}
\end{figure}

\begin{itemize}
    \item in $\dhex{20}$ we set quality weights on nodes and arcs, as done in the rest of the paper: $w(n_{\P})=\varrho(\P)$ and $w(a_{\F})=\varrho(\P_1\cup\P_2)$ for all nodes $n_{\P}$ and arcs $a_{\F}$;
    \item in $\dhex{20c}$ we set constant weights, using METIS in its original mode: $w(n_{\P})=1$ and $w(a_{\F})=1$ for all $n_{\P}$, $a_{\F}$;
\end{itemize}
Different weights produce optimized meshes with similar numbers of vertices and elements, but organized in different configurations.
Using constant weights, we optimize the mesh following the METIS internal criterion for graph partitioning, which essentially tries to agglomerate elements in compact and convex configurations.
This is conceptually not so far from what the quality indicator~\eqref{eq:indicator:3D} suggests, in particular, $\rhov{1}$ and $\rhov{2}$ promote exactly such types of polytopes.
Therefore, it is reasonable that the output meshes will not be so different from each other, and we want to make sure that the use of the quality indicator has an impact.

We see how $\dhex{20}$ produces smaller $\errL$ values (and $\errH$ values are analogous) than $\dhex{20c}$, and this difference remains constant when increasing $k$.
Moreover, meshes in $\dhex{20}$ contain fewer DOFs than those in $\dhex{20c}$, as indicated by the position of the dots in the plot.
It is therefore important to set the weights using a quality criterion, also considering that the weights computation time is negligible compared to the time needed for METIS to partition the graph.

\section{CAD Application}
\label{sec:cad}
\begin{figure}[!h]
  \centering
  \begin{subfigure}[]{.48\textwidth}
      \centering
    \includegraphics[width=\linewidth]{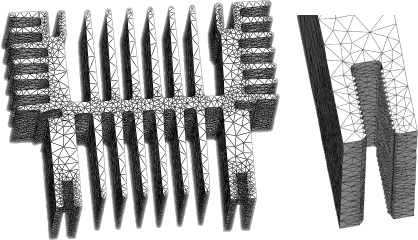}
      
      \vspace{.3cm}
     \includegraphics[width=.8\linewidth]{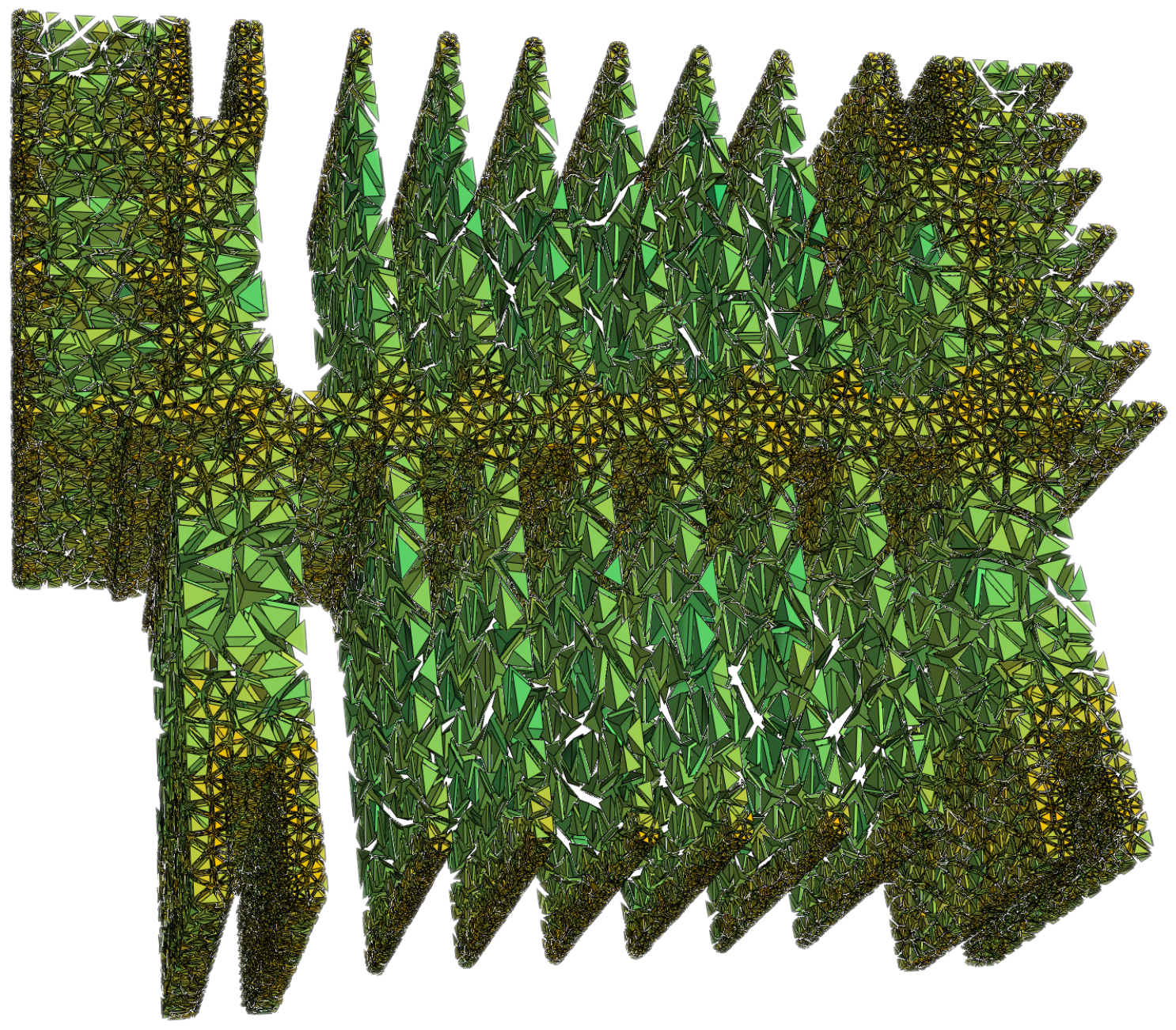}
      \caption{$\dcad{}$, 122882 elements}
      \label{fig:cad100}
  \end{subfigure}
  \hspace{.1cm}\vline\hspace{.1cm}
  \begin{subfigure}[]{.48\textwidth}
      \centering
    \includegraphics[width=\linewidth]{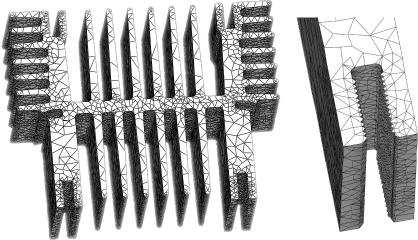}
   
      \vspace{.3cm}
       \includegraphics[width=.8\linewidth]{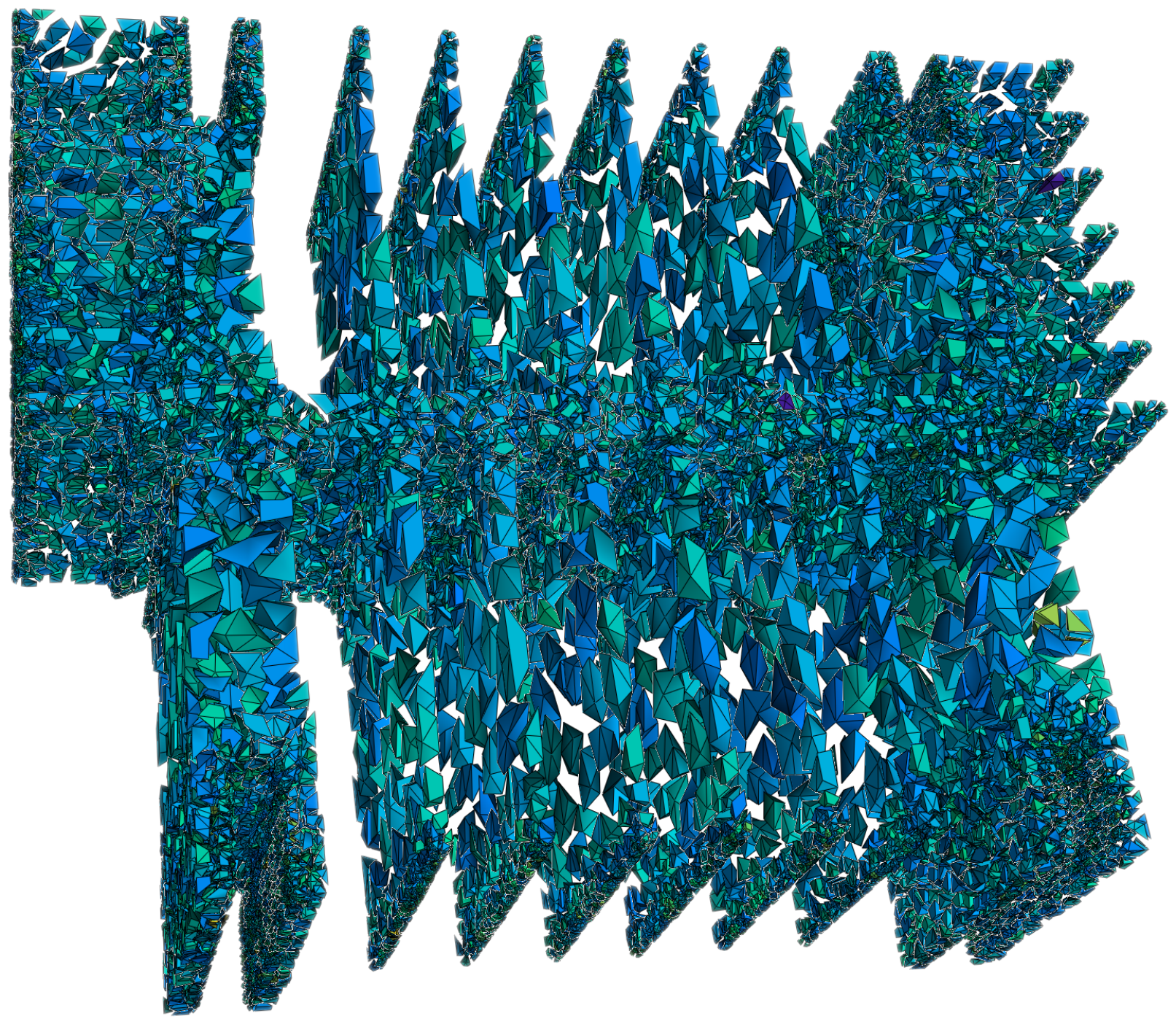}
      \caption{$\dcad{20}$, 24715 elements}
      \label{fig:cad20}
  \end{subfigure}
  \caption{$\dcad{}$ model optimized with $\param=20$, with a close-up on a small-scale detail (top) and elements colored w.r.t. their quality (bottom).}
  \label{fig:cad}
\end{figure}

In this last section we test our algorithm over a real CAD model of a \textit{heat sink}, downloaded from the repository \textit{Traceparts}\footnote{\tiny{\url{www.traceparts.com/it/product/fischer-elektronik-gmbh-co-kg-strangkuhlkorper-fur-einrasttransistorhaltefeder?CatalogPath=TRACEPARTS\%3ATP014002\&Product=34-23062017-133913\&PartNumber=SK\%20593\%2025\&corid=0737fe6b-01e1-4714-6dec-a98692f099d7}}}.
Using \textit{fTetWild}~\cite{hu2020fast}, we re-meshed the original triangulation obtaining a higher quality surface mesh, then we generated a tetrahedralization of the inside.
The considered model is particularly complex because it presents three different levels of detail, which determine the generation of elements with variable size.
The resulting mesh, noted $\dcad{}$, is visible in Figure~\ref{fig:cad}, together with its optimization with parameter $\param=20$.
In the ``exploded'' version, elements are colored with respect to their quality, from yellow ($\varrho_3=1$) to blue ($\varrho_3=0$).
While $\dcad{}$ contains only tets, which are very high quality, in $\dcad{20}$ the $80\%$ of the elements have been agglomerated into generic polyhedra.
The local quality of $\dcad{20}$ is therefore generally lower, and this is an unavoidable side-effect when optimizing high-quality meshes.
However, our guess is that the loss in local quality is negligible when compared to the gain in terms of DOFs.

\subsection{Time-Dependent Problem}
\label{subsec:cad:results}
\begin{figure}[!h]
  \centering
  \begin{subfigure}[]{\textwidth}
      \centering
      \includegraphics[width=.3\textwidth]{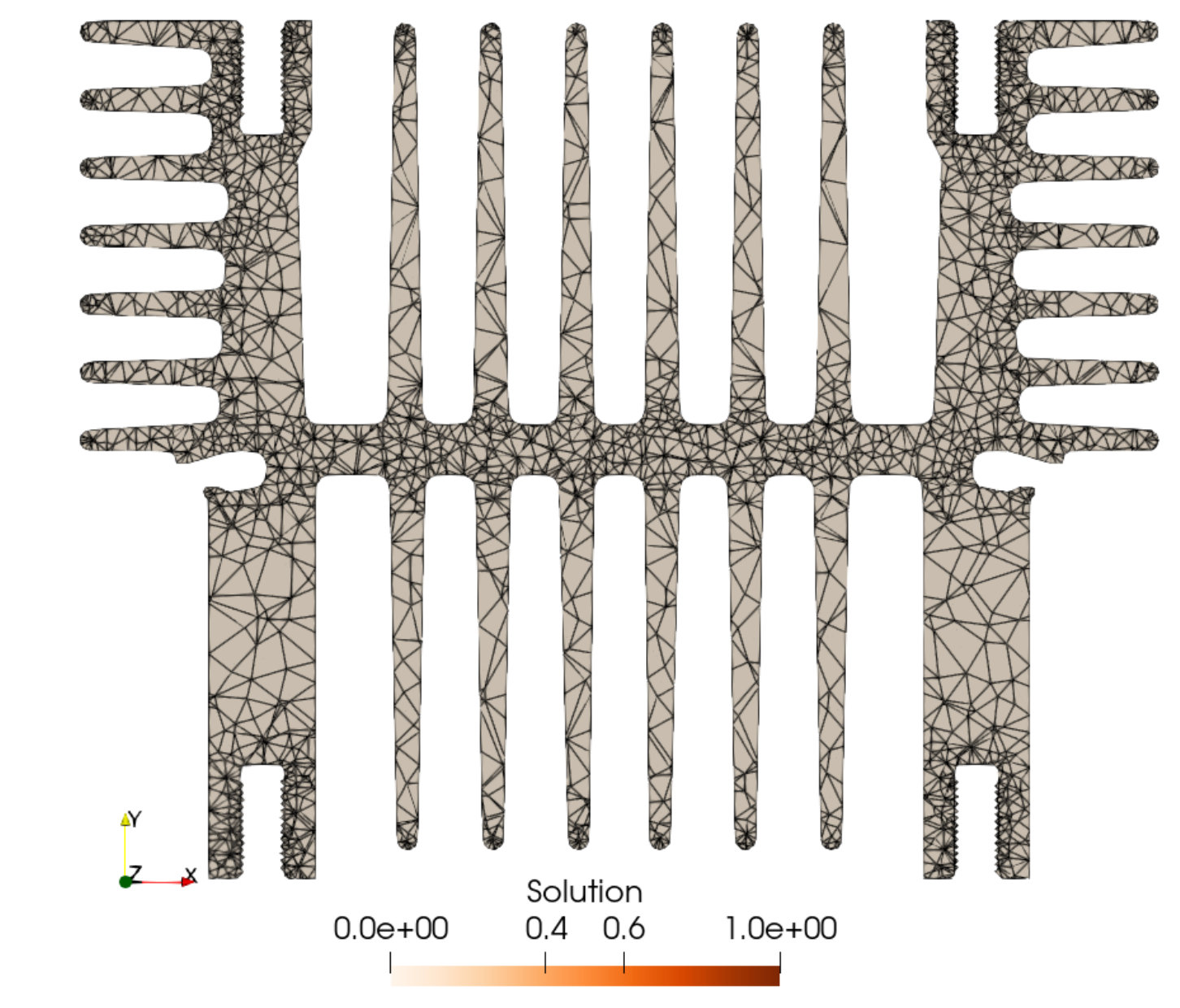}
      \includegraphics[width=.3\textwidth]{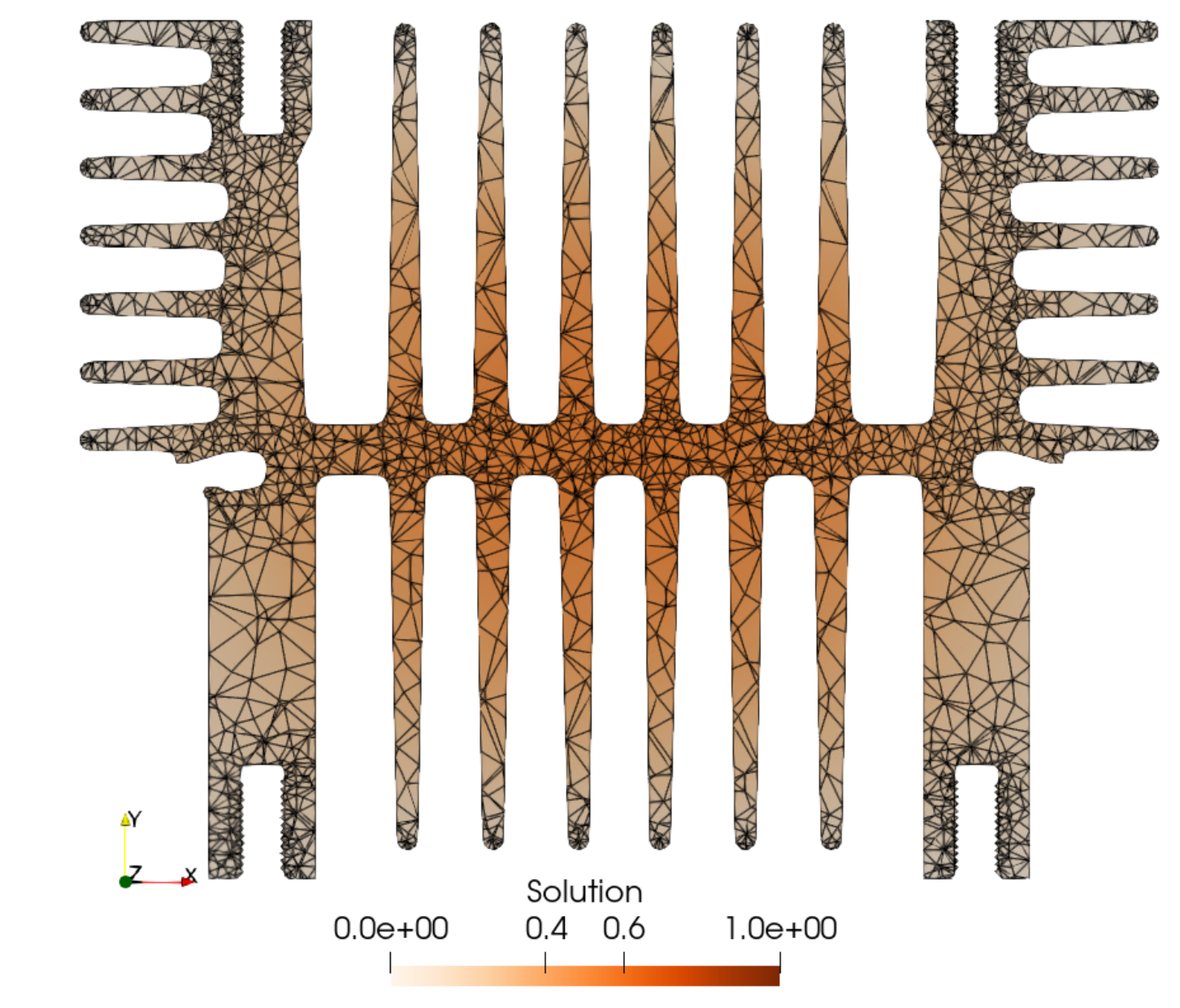}
      \includegraphics[width=.3\textwidth]{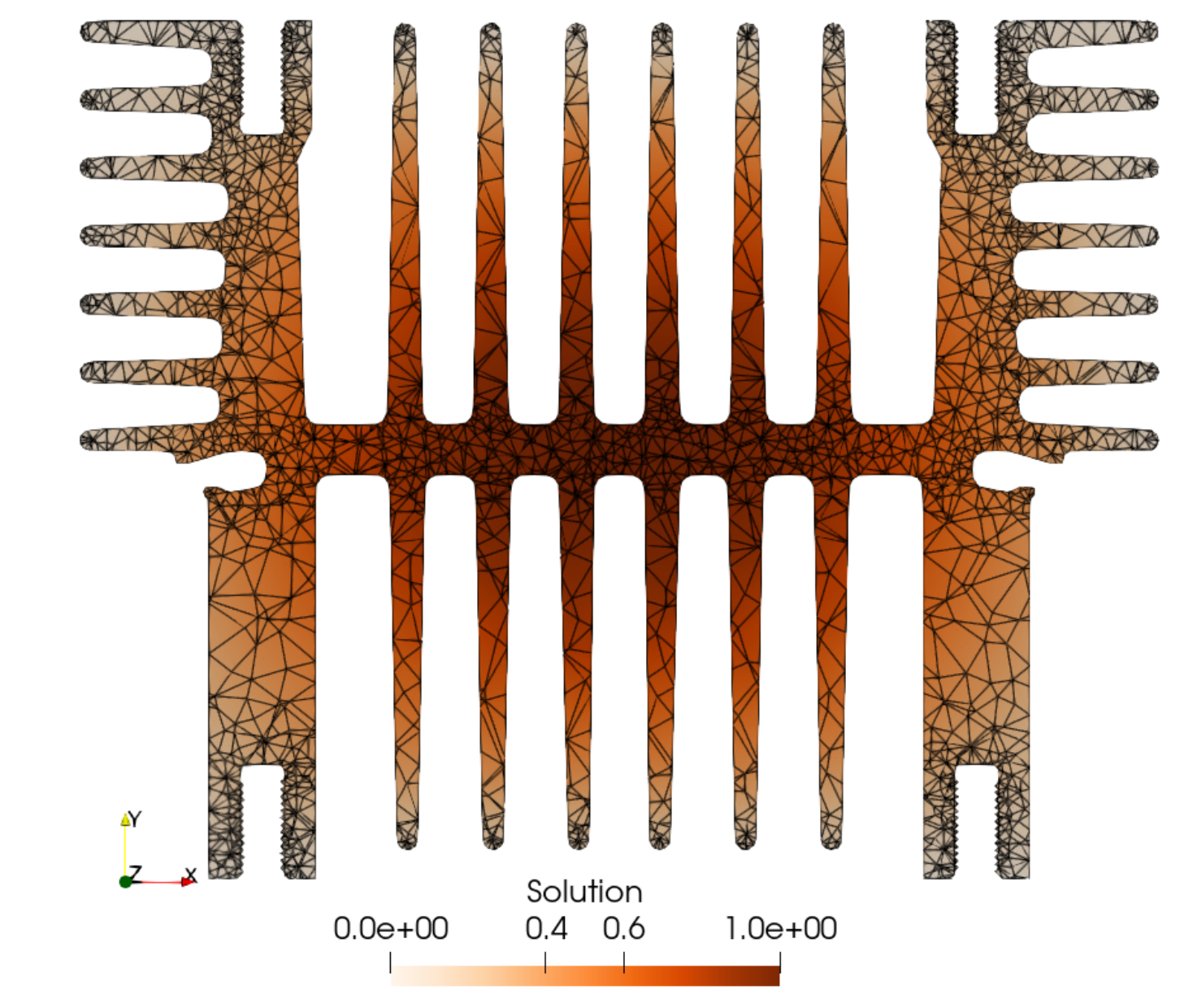}
      \caption{$\dcad{}$}
      \label{fig:cad:simulation100}
  \end{subfigure}
  \vspace{.1cm}
  
  \begin{subfigure}[]{\textwidth}
      \centering
      \includegraphics[width=.3\textwidth]{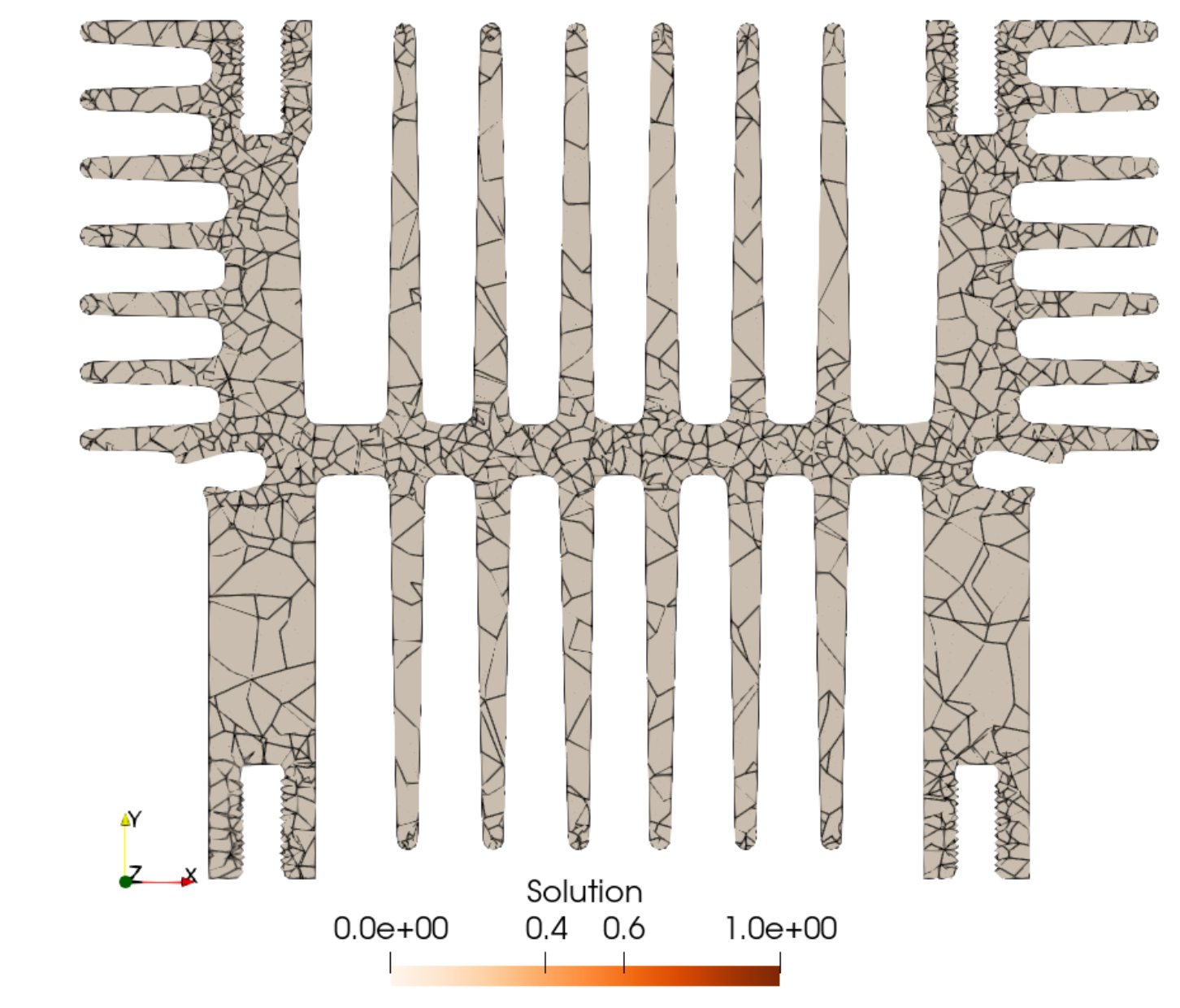}
      \includegraphics[width=.3\textwidth]{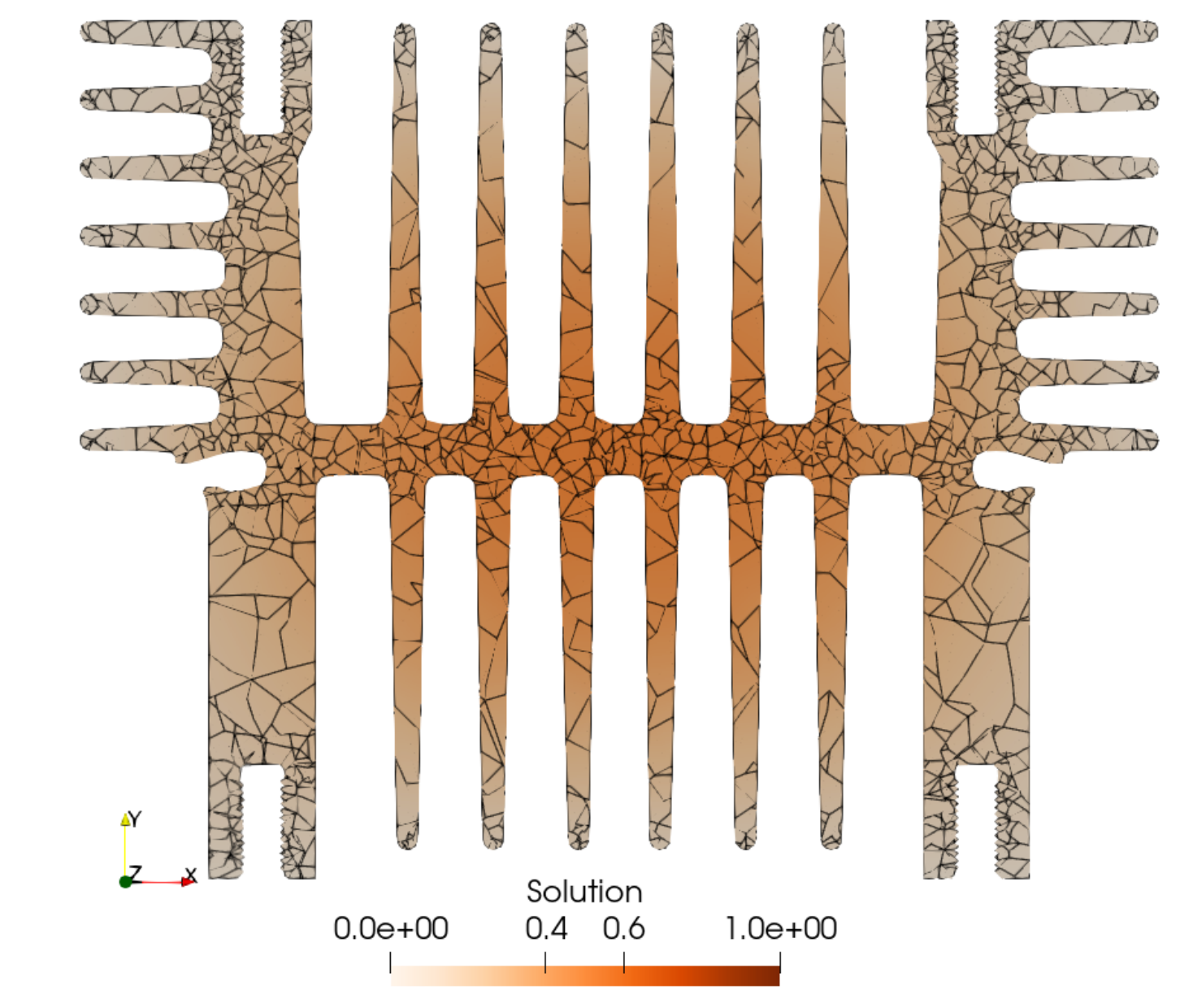}
      \includegraphics[width=.3\textwidth]{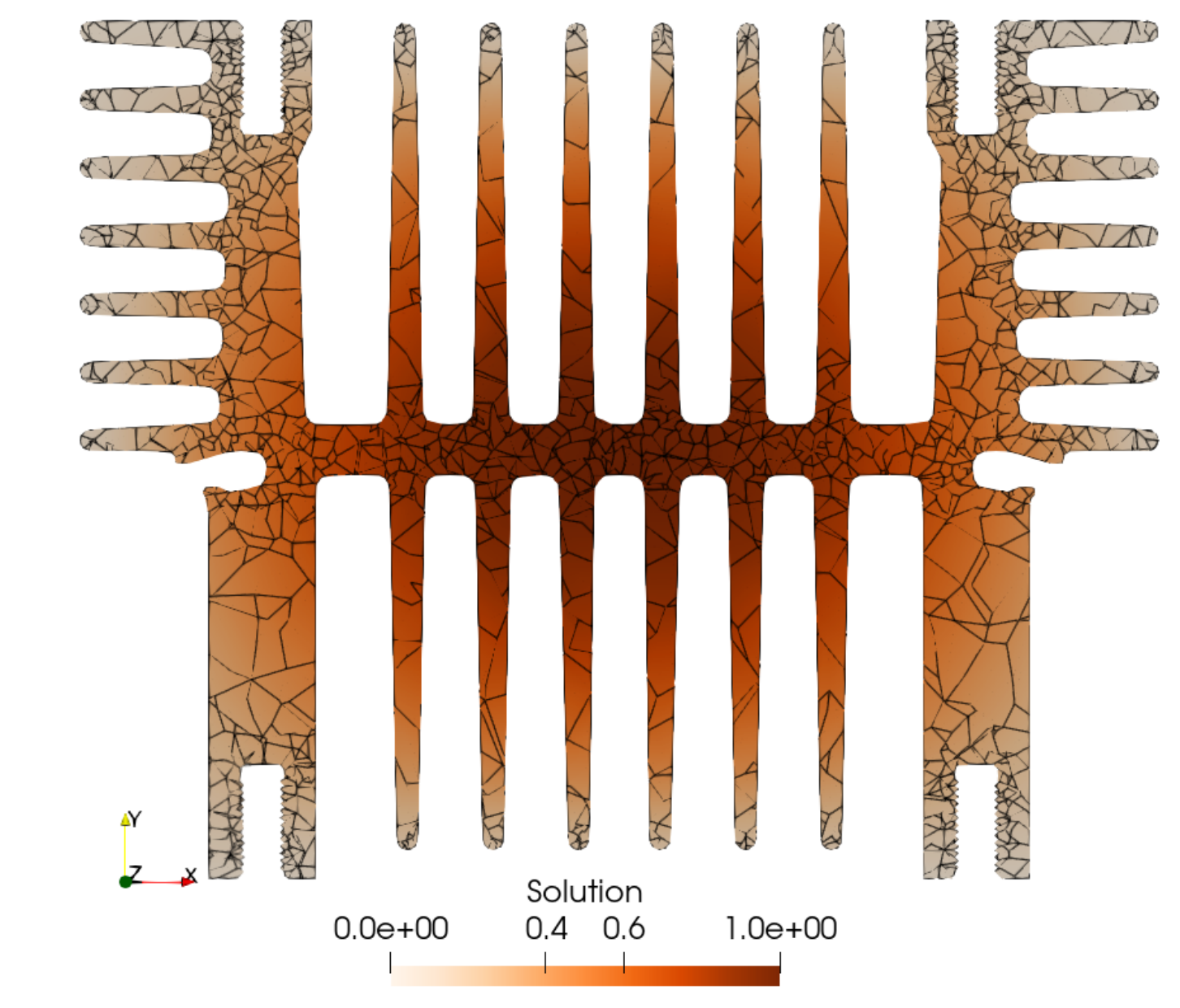}    
      \caption{$\dcad{20}$}
      \label{fig:cad:simulation20}
  \end{subfigure}
  \caption{Simulation with $k=1$ over models $\dcad{}$ (top) and $\dcad{20}$ (bottom). From left to right, $t=0$, $t=0.5$, $t=1$.}
  \label{fig:cad:simulation}
\end{figure}
\begin{figure}[!h]
  \centering
  \includegraphics[width=.7\textwidth]{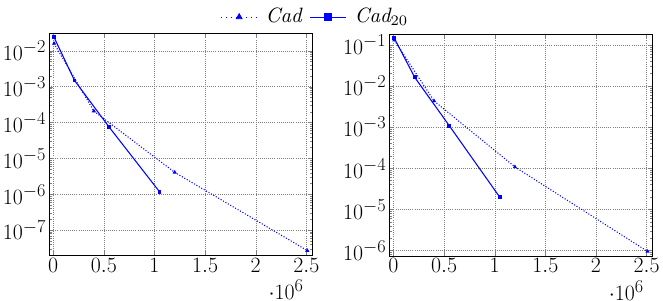}
  \caption{Performance of the VEM with $k=1,2,3,4$ at $t=1$ on models $\dcad{}$ and $\dcad{20}$. We measure $\errL$ (left) and $\errH$ (right) with respect to the number of DOFs in \textit{semilog-$y$} scale.}
  \label{fig:cad:convergence}
\end{figure}
\begin{table}[!h]
  \centering
  \caption{Differences between $\dcad{}$ and $\dcad{20}$ for $k=1,2,3,4$ at $t=1$.}
  \label{tab:cad:delta}
  \includegraphics[width=\textwidth]{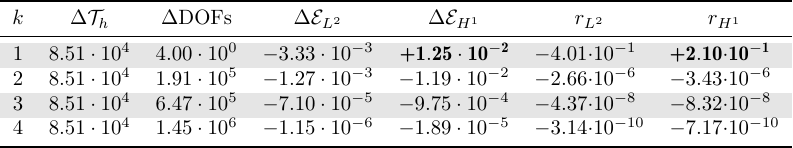}
\end{table}

Differently from Section~\ref{sec:convergence}, we now have a single mesh and we measure the advantages of solving multiple problems over a smaller and better (in the sense of quality) mesh.
We solve the following time-dependent problem, defined on $\dcad{}$ and its optimized version $\dcad{20}$, for a number of time values between 0 and 1:
\begin{eqnarray*}
    \frac{\partial u}{\partial t} - \Delta u &= f &\text{in } \Omega,\ \forall t \in [0, 1],\\
    u(t, x) &= 0 &\text{in } \partial \Omega,\ \forall t \in [0, 1],\\
    u(0, x) &= 0 &\text{in } \bar{\Omega},
\end{eqnarray*} 
The two meshes share the same bounding box with centroid $\xv$, origin $(x_0,y_0,z_0)\cong(-0.035,0.0,-25.002)$, and size $(x_s,y_s,z_s)\cong(75.58,60.01,25.05)$.
Given a constant $\tilde{u}_0$ computed such that $u=1$ when $t=1$, we seek for the exact solution:
\begin{equation*}
    u(\xv, t) = t \tilde{u}_0 (x - x_0)(y - y_0)(z - z_0)(x-(x_0 + x_s))(y-(y_0 + y_s))(z-(z_0 + z_s)).
\end{equation*}
We use a backward Euler scheme with 10 time steps; visual results are presented in Figure~\ref{fig:cad:simulation}.

In Figure~\ref{fig:cad:convergence} we plot $\errL$ and $\errH$ for $k=1,2,3,4$ relatively to $t=1$, but these errors remain very similar for all time steps.
Instead of refining the mesh at each iteration, as done in Section~\ref{sec:convergence}, we increase the order of the VEM scheme.
The original meshes produce more accurate results point-wise, i.e., the $i$-th dot of $\dcad{}$ is higher than the $i$-th dot of $\dcad{20}$ in the plot.
However, the difference in error is much smaller than the difference in terms of DOFs.
Both plots converge linearly in \textit{semilog-y} scale as expected, but the line relative to $\dcad{20}$ is significantly below the $\dcad{}$ one.

In Table~\ref{tab:cad:delta} we can see how $\delD$ grows with $k$ while $\delL$ and $\delH$ decrease, all with an exponential rate.
Therefore, despite $\delL$ and $\delH$ being negative, their influence is smaller and smaller.
The last column represents the ratio $r_{\LTWO}=\delL/\delD$, and analogously for $r_{\HONE}$.
When we compare two meshes in the plot we have an $x$-shift represented by $\delD$, and an $y$-shift represented by $\delL$ (or $\delH$).
The ratio $r=dx/dy$ then represents the slope of the line connecting the two meshes.
The smaller $|r|$, the more is convenient to pass from one mesh to the other, because the difference in DOFs overcomes the difference in error.
For instance, for $k=3$ we have $|r_{\LTWO}|\sim~10^{-8}$: we are losing a factor $10^{-8}$ in accuracy for every DOF that we removed from the original mesh.

\subsection{Computational Time}
\label{subsec:cad:time}
\begin{table}[!h]
  \centering
  \caption{Differences between computational times, in seconds, for the time-dependent problem with $k=3$ solved on $\dcad{}$ and $\dcad{20}$.
  In the last column, percentage of the total time required by the optimization of the mesh.}\label{tab:cad:time}
  \includegraphics[width=.9\textwidth]{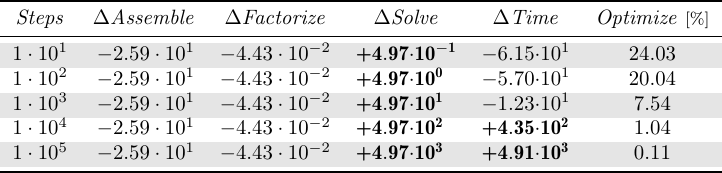}
\end{table}

In Table~\ref{tab:cad:time} we report a comparison of the computational times required for the main operations involved in the process, for the case $k=3$.
As in the rest of the paper, the comparison is made subtracting the time relative to $\dcad{20}$ from the time relative to $\dcad{}$, hence a positive $\Delta$ means a gain in time when using the optimized mesh.
Computations have been performed using a 3.60GHz Intel Core i7 processor with 16 cores and 32GB of RAM.

The first column contains the number of time steps considered, from 10 to $10^6$.
The times for importing the mesh structure in the solver, computing the geometric properties, and assembling the linear system, are aggregated in \textit{Assemble}, while \textit{Factorize} indicates the time for factorizing the global matrix of the problem.
These operations are performed only once, therefore these times are independent of the number of time steps.
The differences $\Delta$\textit{Assemble} and $\Delta$\textit{Factorize} are negative because elements in optimized meshes generally have more complicated shapes than the simple tets of the original ones, despite being fewer.
We indicate by \textit{Solve} the time required for the solution of the problem, resulting from the sum of \textit{Steps} iterations.
The solution time is faster on the optimized meshes because the problem itself becomes much smaller, and this gain gets multiplied by the number of time steps.
\textit{Time} indicates the total time, resulting from the sum of the previous ones and including also the optimization of the mesh.
Note that $\Delta$\textit{Time} cannot be expressed in terms of \eqref{eq:delta_t}, as each $\delT$ already contains the time for optimizing the mesh. 
This $\Delta$ is initially negative, but it grows with the number of time steps because the initial fixed costs (optimizing the mesh, assembling the system, factorizing the matrix) become negligible with respect to \textit{Solve}.
In particular, we report the impact of the optimization on the total time in the last column.

We could observe a difference in the \textit{Factorize} cost with respect to $\dtet{}$ and $\dhex{}$ (Table~\ref{tab:results} and Table~\ref{app:results:3D}).
In those cases, the factorization step is significantly cheaper on the optimized mesh, resulting in a positive $\delT$ for $k>1$ even if we were solving a single iteration of the problem.
We associate this fact with the quality of the original mesh: $\dtet{}$ and $\dhex{}$ contain very bad-shaped elements, while $\dcad{}$ is made of good-quality tets.
If the starting mesh is already good, there is less space for improvements, and our optimization becomes effective only after a certain number of time steps.

In particular, there exists a critical number of time steps $\bar{t}$ such that for less than $\bar{t}$ steps, computations are faster on the original mesh, while for more than $\bar{t}$ steps the optimized mesh becomes advantageous.
This number decreases as $k$ grows: we have $\bar{t}=2.63\cdot10^4$ for $k=2$, $\bar{t}=1.25\cdot10^3$ for $k=3$, and $\bar{t}=7.63\cdot10^2$ for $k=4$.
For $k=1$ times are essentially equivalent, thus we could not identify a single $\bar{t}$ value.
However, the presented calculations were made using a code written and optimized for meshes with convex elements, therefore we believe further improvements can be made with solvers optimized for generic non-convex cells.

\section{Conclusions}
\label{sec:conclusions}
We presented an algorithm for optimizing the size of a mesh with respect to its quality, in the context of VEM simulations.
The algorithm is able to significantly reduce the total number of mesh elements and the number of DOFs, in particular in high-order formulations, while preserving the VEM optimal convergence rates. 
These effects lead to a notable decrease in the computational effort required for obtaining the numerical solution.
A scenario in which this decrease becomes particularly interesting is the class of time-dependent problems: in such cases, we optimize the mesh only once and then save time in each iteration.
We also observed how, if the input mesh has a particularly low quality, the optimization can locally improve the shape of the elements.
This improvement allows to control the stiffness matrix condition number and therefore to recover the optimal VEM convergence rate.
The mesh optimization algorithm is therefore a powerful tool to achieve faster and more stable VEM simulations.

One potential limitation of the algorithm is that it only evaluates elements pairwise. 
When considering an element $\P$ with neighboring elements $\P'$ and $\P''$, the algorithm separately assesses the potential quality of $\P\cup\P'$ and $\P\cup\P''$, without considering the overall quality of $\P\cup\P'\cup\P''$. 
Consequently, some small mesh elements may still persist in the optimized final mesh, particularly around small-scale shape features of the boundary.
This issue can be controlled by subdividing the optimization into smaller steps (i.e., running the algorithm twice with $\param=40$ and $\param=25$ instead of once with $\param=10$).
However, the numerical findings presented in this study suggest that such action is generally not necessary for obtaining satisfactory results.

We believe that an interesting aspect of our algorithm is its generality, as it can be applied to any type of mesh.
It is also modular: its key ingredients, the mesh quality indicator and the graph partitioning algorithm, may be easily replaced or adjusted for adapting the optimization to specific needs.
With a different quality indicator, we could optimize the mesh to be used with numerical methods other than the VEM or include in the quality criterion features of the numerical problem (e.g., anisotropies).
We refer the reader to \cite{sorgente2023survey} for more quality indicators.
Other algorithms for graph partitioning may be more efficient, see the parallel version of METIS, ParMETIS~\cite{karypis1997parmetis}.
In some situations instead, it may be useful to have a partitioning method that finds a local maximum of the global quality of the mesh and automatically decides the optimal number of partitions, as it is done in \cite{sorgente2023mesh} using GraphCut~\cite{boykov2001fast}.
Future research will investigate these aspects more deeply, as well as try to adapt a similar approach to mesh refinement strategies.

\section*{Acknowledgments}
This work was carried out within the framework of the project ``RAISE - Robotics and AI for Socio-economic Empowerment'' and has been partially supported by European Union - NextGenerationEU. 
However, the views and opinions expressed are those of the authors alone and do not necessarily reflect those of the European Union or the European Commission. 
Neither the European Union nor the European Commission can be held responsible for them.

F. Vicini, S. Berrone, and G. Manzini are members of the Gruppo Nazionale Calcolo Scientifico-Istituto Nazionale di Alta Matematica (GNCS-INdAM); 
this work has been partially supported by INdAM-GNCS Project CUP: E53C23001670001, 
by the MUR PRIN project 20204LN5N5\_003,
and by the European Union through the project Next Generation EU, PRIN 2022 PNRR project P2022BH5CB\_001 CUP:E53D23017950001.

\bibliographystyle{unsrt}
\bibliography{biblio}

\appendix
\section{Supplementary Material}
\label{app}
We collect in Tables~\ref{app:dataset:optimization:3D}-\ref{app:dataset:optimization:2D} the analogues of Table~\ref{tab:optimization} for datasets $\dhex{}$, $\dtri{}$, $\dqua{}$, and their optimizations with $\param=40$ and $\param=20$.
In Tables~\ref{app:results:3D}-\ref{app:results:2D} we report the analogues of Table~\ref{tab:results} for datasets $\dhex{20}$, $\dtri{20}$, and $\dqua{20}$.

\begin{table}[htbp]
    \centering
    \caption{Analogue of Table~\ref{tab:optimization} for dataset $\dhex{}$.
    Number of vertices, faces and elements of each mesh; number of internal DOFs for $k=1,2,3$; mesh quality~\eqref{eq:indicator:mesh}.}
    \label{app:dataset:optimization:3D}
    \includegraphics[width=.9\textwidth]{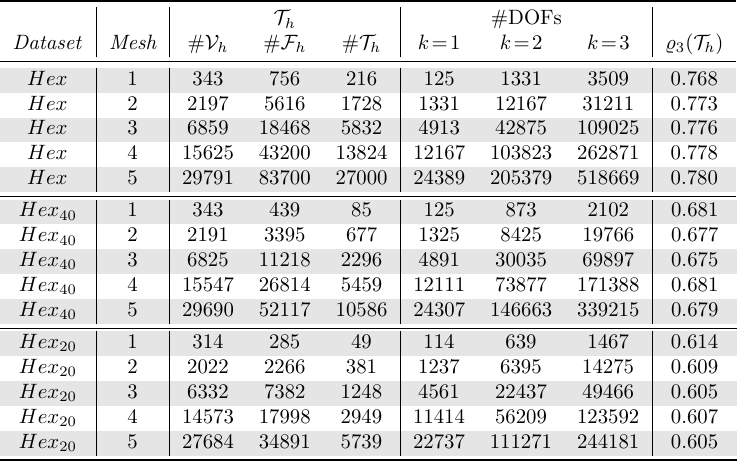}
\end{table}

\begin{table}[htbp]
    \centering
    \caption{Analogue of Table~\ref{tab:optimization} for datasets $\dtri{}$ and $\dqua{}$.
    Number of vertices, edges and elements of each mesh; number of internal DOFs for $k=1,2,3$; mesh quality~\eqref{eq:indicator:mesh}.}
    \label{app:dataset:optimization:2D}
    \includegraphics[width=.9\textwidth]{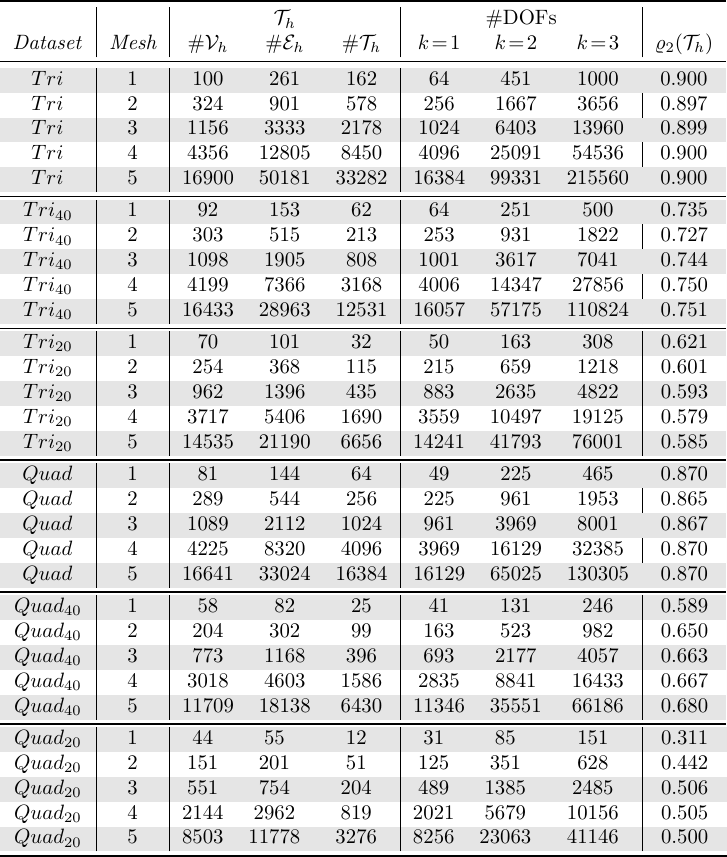}
\end{table}

\begin{table}[htbp]
    \centering
    \caption{Analogue of Table~\ref{tab:results} for dataset $\dhex{}$.
    Difference between the number of DOFs~\eqref{eq:delta_dofs}, the $\errL$ error~\eqref{eq:delta_L2}, the $\errH$ error~\eqref{eq:delta_H1}, and the computational time~\eqref{eq:delta_t} w.r.t. the original mesh.} 
    \label{app:results:3D}
    \includegraphics[width=.9\textwidth]{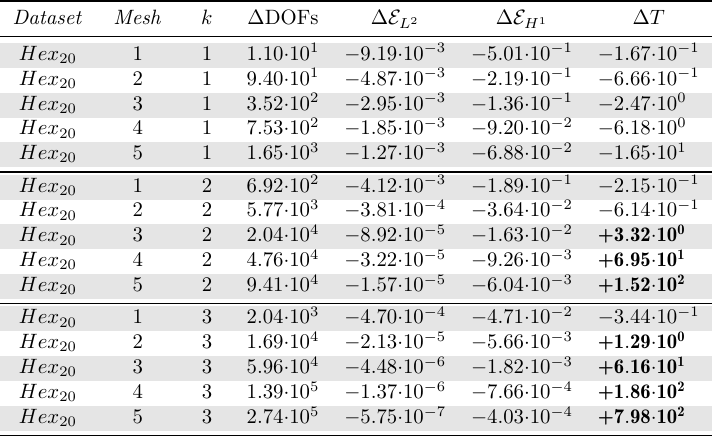}
\end{table}

\begin{table}[htbp]
    \centering
    \caption{Analogue of Table~\ref{tab:results} for datasets $\dtri{}$ and $\dqua{}$.
    Difference between the number of DOFs~\eqref{eq:delta_dofs}, the $\errL$ error~\eqref{eq:delta_L2}, the $\errH$ error~\eqref{eq:delta_H1}, and the computational time~\eqref{eq:delta_t} w.r.t. the original mesh.} 
    \label{app:results:2D}
    \includegraphics[width=.9\textwidth]{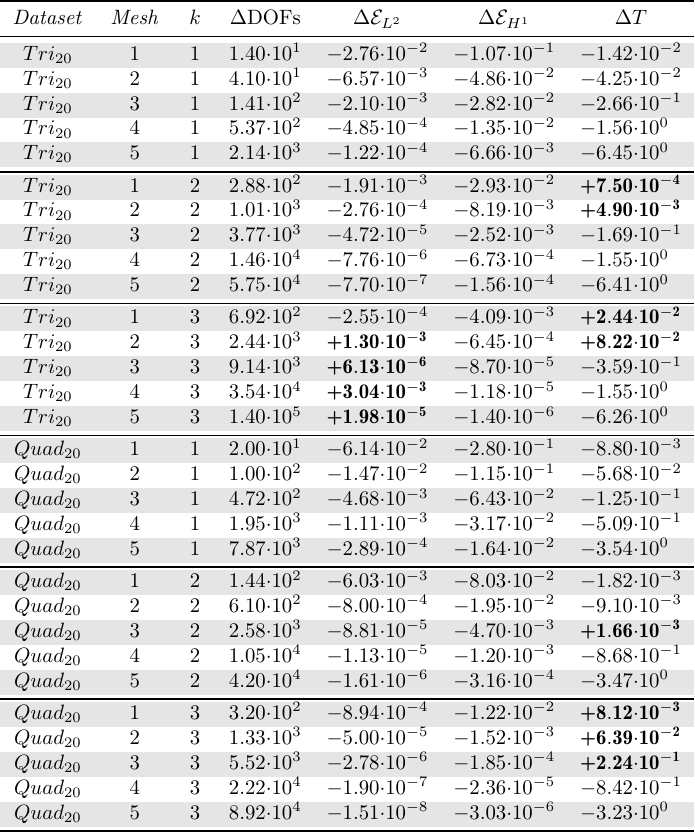}
\end{table}

\end{document}